\documentclass[leqno]{amsart}
\usepackage{amsfonts,amssymb,amsmath,amsgen,amsthm}

\theoremstyle{plain}
\newtheorem{theorem}{Theorem}[section]

\newtheorem{assumption}[theorem]{Assumption}
\newtheorem{lem}[theorem]{Lemma}

\newtheorem{prop}[theorem]{Proposition}

\newtheorem{rem}[theorem]{Remark}

\newtheorem{example}[theorem]{Example}

\def\C{{\mathbf C}}% complex numbers
\def\R{{\mathbf R}}% real numbers
\def\N{{\mathbf N}}% nonnegative integers
\def\({\left(}
\def\){\right)}
\def\<{\left\langle}
\def\>{\right\rangle}

\def\d{{\partial}}

\numberwithin{equation}{section}

\def\C{{\bf C}}

\def\N{{\bf N}}
\def\R{{\bf R}}

\def\op_#1{\mathrel{\mathop{{\rm op}_{#1}}}}

\def\build#1_#2^#3{\mathrel{
\mathop{\kern 0pt#1}\limits_{#2}^{#3}}}

\def\td_#1,#2{\mathrel{
\mathop{\build\longrightarrow_{#1\rightarrow #2}^{}}}}

\def\lim_#1,#2{\mathrel{
\mathop{\build{\rm lim}_{#1\rightarrow#2}^{}}}}

\def\limsup_#1,#2{\mathrel{
\mathop{\build{\rm limsup}_{#1\rightarrow#2}^{}}}}

\def\liminf_#1,#2{\mathrel{
\mathop{\build{\rm liminf}_{#1\rightarrow#2}^{}}}}

\def\aref#1{(\ref{#1})}

\newtheorem{defi}{Definition}

\def\1{{\bf 1}}
\def\0{{\bf 0}}

\def\ni{\noindent}
\def\op{{\rm op}}

\def\d{{\rm d}}

\begin{document}

\title[Analysis of the Energy decay of a  Degen. Thermoelasticity System]{Analysis of the Energy Decay of  a Degenerated Thermoelasticity System}
\author[A. Atallah]{Amel Atallah-Baraket}
\address[A. Atallah]{Dpt. de Math\'ematiques, Facult\'e des Sciences de Tunis, Universit\'e de Tunis El Manar, 1060, Tunis, Tunisie.}
\email{amel.atallah@fst.rnu.tn}
\author[C. Fermanian]{Clotilde~Fermanian Kammerer}
\address[C. Fermanian]{LAMA UMR CNRS 8050,
Universit\'e Paris Est - Cr\'eteil\\
61, avenue du G\'en\'eral de Gaulle\\
94010 Cr\'eteil Cedex\\ France.}
\email{Clotilde.Fermanian@u-pec.fr}
\date{}
\maketitle

{\bf Abtract}: {\it In this paper, we study a system of thermoelasticity with a degenerated second order operator in the Heat equation. We analyze the evolution of the energy density of a family of solutions. We consider two cases: when the  set of points  where the ellipticity of the Heat operator fails is  included in a hypersurface and when it is an open set. 
In the first case and under special assumptions, we prove that  the evolution of the energy density is the one of  a damped wave equation: propagation along the rays of geometric optic and damping according to a microlocal process. In the second case, we show that  the energy density propagates along rays which are  distortions of the rays of geometric optic.}

\section{Introduction}

We consider $\Omega$ an open subset of $\R^d$ and  the following system of thermoelasticity
\begin{eqnarray}\label{eq:thermo}
 \left\{\begin{array}{lll}
\partial^2_t u - \Delta u + \nabla\cdot  (\gamma(x) \theta) = 0, \;\; (t,x) \in\R^+ \times \Omega,\\
\partial_t \theta - \nabla. (B(x) \nabla \theta) + \gamma(x). \nabla
\partial_tu = 0, \;\;  (t,x) \in\R^+ \times \Omega,\\
u_{|t = 0} = u_0\;\;{\rm and }\;\;
\partial_tu_{|t=0} = u_1 \; \mbox{ in} \;\Omega\\
\theta_{|t=0} = \theta_0 \; \mbox{ in} \;\Omega,\\
u_{|\partial\Omega}=0,\;\;\partial_tu_{|\partial\Omega}=0\;\;{\rm and}\;\;\theta_{|\partial\Omega} =0.
\end{array}\right.
\end{eqnarray} where $u$ and $\theta$ are scalar real-valued
functions. The matrix-valued function $x\mapsto B(x)$   and the vector-valued function
$x\mapsto\gamma(x)$ are
supposed to be defined on $\R^d$  (and thus on $\overline\Omega$) and to depend smoothly on the variable $x\in\R^d$.
 The matrix $B(x)$ is  assumed to be symmetric, non-negative:
there exists $C_2>0$ such that
\begin{eqnarray}\label{eq:B(x)}
\forall\; (x,\xi) \in \overline \Omega\times \R^{d},\;\; 0\leq B(x) \xi.\xi \leq C_2
|\xi|^2.
\end{eqnarray}
Note that we may have ${\rm det} \,B(x)=0$.

$ $

\ni  System \aref{eq:thermo} has an energy
$$E(u,\theta,t) := \int_\Omega |\partial_t u(t,x)|^2 dx + \int_\Omega \left|
 \nabla u(t,x)\right|^2 dx +\int_\Omega\theta^2(t,x) dx$$ which decreases in
time according to
\begin{eqnarray}\label{eq:energy}
E(u,\theta,t)  - E(u,\theta,0) = -2 \int^t_0\int_\Omega B(x) \nabla \theta(s,x). \nabla
\theta(s,x) dsdx \leq 0.
\end{eqnarray}

\ni Equation  \aref{eq:energy} gives {\it a priori}
estimates for initial data $u_0 \in H^1_0(\Omega) , u_1 \in L^2(\Omega)$
and $ \theta_0 \in L^2(\Omega) $ and yields by classical arguments
the existence of a unique solution $(u,\theta) \in \mathcal{C}^0(
\R^+, H^1_0(\Omega))
 \cap \mathcal{C}^1( \R^+, L^2(\Omega)) \times \mathcal{C}^0( \R^+, L^2(\Omega)).$
We are interested in characterizing the way the energy decay: our aim is to describe the weak limits of the energy densities
$$e_n(t,x)=|\partial_t u_n(t,x)|^2 + \left|\nabla_x
u_n(t,x)\right|^2 +(\theta_n(t,x) )^2,\;\;n\in\N,$$
associated with families of solutions $(u_n,\theta_n)_{n\in\N}$ of \aref{eq:thermo} corresponding to families of initial data $(u_{0,n})_{n\in\N}$ in one hand and $(u_{1,n})_{n\in\N}$, $(\theta_{0,n})_{n\in\N}$ in the other hand, uniformly bounded in  $H^1(\Omega)$ and $L^2(\Omega)$ respectively. Without loss of
generality, we suppose that $(u_1^n)_{n\in \N}$ and $(\theta_0^n)_{n\in\N}$ goes  to
zero weakly in $L^2(\Omega)$ and that $(u_0^n)_{n\in\N}$ goes  to zero weakly in
$H^1(\Omega)$.

$ $

Main results on the subject are devoted to the situation where the matrix $B$ is positive.
It is  known since the work of Dafermos \cite{D} that the energy decays and  the description  of this decay has been the subject of several contributions. In particular,   in~\cite{LZ},    G.~Lebeau and E.~Zuazua \cite{LZ} have proved that the decay of solutions to \aref{eq:thermo} is not uniform.
They crucially use a result of Henry, Perissinoto and Lopes~\cite{HPL}  which show that the semigroup associated to \aref{eq:thermo} is equal up to a compact operator to the semi-group of a system consisting of a damped wave equation coupled with the equation on the temperature $\theta$:
\begin{equation}\label{eq:dw}
\left\{\begin{array}{l}
\partial^2_t u - \Delta  u + \Gamma u= 0, \;\; (t,x) \in\R^+ \times \Omega\\
\partial_t \theta - \nabla \cdot (B(x) \nabla \theta) + \gamma(x). \nabla
\partial_tu = 0.
\end{array}\right.
\end{equation}
The operator $\Gamma $ is a damping operator given by
$$\Gamma= G^* Q^{-1} G,\;\;G=\gamma(x)\cdot \nabla \;\;{\rm and}\;\; Q=\nabla\cdot  \left(B(x)\nabla \cdot \right).$$
The behaviour of the energy density $|u_n(t,x)|^2$
associated with families of solutions  of this damped wave equation (for initial  data $(u_{0,n})_{n\in\N}$ and $(u_{1,n})_{n\in\N}$ uniformly bounded in  $H^1(\Omega)$ and $L^2(\Omega)$ respectively) can be studied in the same manner than in the paper of  G.~Lebeau and N.~Burq (see \cite{L1}, \cite{L2} and \cite{BL}).   In these papers, the damping is somehow different than the one of~\aref{eq:dw}, however their method applies. One obtains that the energy propagates along the  rays of geometric optic associated with the wave operator $\partial_t^2-\Delta$ with a damping depending simultaneously on the position and the speed of the trajectory. These trajectories are trajectories of the phase space $\left(x(t),\omega(t)\right)$ and the damping is microlocal: it is given by the principal symbol of the operator $\Gamma$,  the function $(x,\omega)\mapsto-{(\gamma(x)\cdot\omega)^2\over B(x)\omega\cdot\omega}$.

$ $

The method of~\cite{BL} is based on the use of microlocal defect measures in the spirit of the article \cite{GL} by P.~G\'erard and E.~Leichtnam who are at the origin of that method.
Then, similar works have  been achieved for the Lam\'e system in \cite{BL}, for the equations of magnetoelasticity in \cite{D1} (see also \cite{D2}) and for the equation of viscoelastic waves by the authors \cite{AF}.

$ $

Of course, the strategy that we have just described, fails if the kernel of $B$  does not reduce to zero and we are interested in this situation: we assume that the set
$$\Lambda=\{(x,\omega)\in\Omega\times {\bf S}^{d-1},\;\;B(x)\omega=0\}$$
is not empty. Then, the operator $\nabla\cdot \left(B(x)\nabla\cdot\right)$ is no longer elliptic  and something else has to be done.
The method we use  to treat the coupling between the temperature~$\theta$ and the amplitude~$u$  is mainly inspired by the analysis of semiclassical systems performed in \cite{GMMP}.
Of course, we recover the result sketched in the previous paragraph when $\Lambda=\emptyset$ and we are also able to extend it to situations where  $\Lambda\not=\emptyset$ provided a {\it weak degeneracy assumption} stated below (see Assumption~\ref{assumption}). This assumption consists first in a geometric assumption: the projection of $\Lambda$ on $\Omega$ is included in a hypersurface $\Sigma$ and  for $(x,\omega)\in\Lambda$, the vector $\omega$ is transverse to $\Sigma$ at the point $x$.  Then, Assumption~\ref{assumption} contains a compatibility relation between the vector $\gamma(x)$ and the matrix $B(x)$: $\gamma(x)\in{\rm Ran } \, B(x)$. With these assumptions, we are able  to prove that the energy density is still damped along the rays of geometric optic even though they pass through $\Sigma$. At the opposite,  if $B(x)=0$ in an open subset $\widetilde\Omega$ of $\Omega$, then the damping disappears and we have transport of the energy along rays which are distorsion of the rays observed before.  Precise statements of our results are given in Section~2 and the organization of the paper is discussed at the end of this section.

$ $

\noindent{\it Notations}. We will say that a sequence $(u_n)_{n\in\N}$ is u.b. in the functional space $F$ if the sequence $(u_n)_{n\in\N}$ is a uniformly bounded family of~$F$.
 We denote by $|X|_{\C^{d+2}}$ the hermitian norm of $X=(X_1,\cdots,X_{d+2})\in\C^{d+2}$: $$|X|^2_{\C^{d+2}}=X_1^2+ \cdots +X_{d+2}^2.$$ Similarly, we will use the notation $(X|Y)$ for the hermitian scalar product of $\C^{d+2}$:\\
$\forall X=(X_1,\cdots,X_{d+2})\in\C^{d+2}$, $\forall Y=(Y_1,\cdots,Y_{d+2})\in\C^{d+2},$,
$$(X|Y)_{\C^{d+2}}=X_1\overline{Y_1}+\cdots + X_{d+2} \overline{Y}_{d+2}.$$

%%%%%%%%%%%%%%%%%%%%%%%%%%%%%%%%%%%%%%%%%%%%%%%%%%%%%

\section{Main results}

In this section, we present our results which crucially rely on the use of microlocal defect measures that we define in the first subsection: the evolution of the energy density is a corollary of the analysis of the behavior of microlocal defect measures associated with the sequence $(u_n,\theta_n)_{n\in \N}$. The second subsection is devoted to the analysis of properties of the thermoelasticity operator that are important for our purpose. Then, in the third subsection, we mainly consider  the situation where the determinant of $B$ vanishes on points of $\Omega$ which are simultaneously included in a hypersurface and in a compact subset of $\Omega$ (so that $B$ is non negative in a neighborhood of $\partial\Omega$). Finally, in the fourth subsection, we discuss what happens if  $B$ vanishes in an open subset of $\Omega$.

\subsection{Microlocal defect measures}

 Microlocal defect measures allow to treat quadratic quantities like the energy density by taking into account microlocal effects. They describe up to a subsequence the limit of quantities of the form $(a(x,D) f_n,f_n)$ where $a(x,D)$ is a pseudodifferential operator and $(f_n)_{n\in\N}$ a u.b. family of $L^2(\Omega)$ (or, more generally, of $H^s(\Omega)$). Recall that the pseudodifferential operator $a(x,D)$ is characterized by its symbol $a(x,\xi)$ which is a smooth function taken, for example, in the space
 ${\mathcal A}_i^m$ of symbols of order $m$: this set contains the functions
  $a=a(x,\xi)$  of ${\mathcal C}^\infty (\Omega\times \R^{d})$ such that
$a$ is compactly supported in $\Omega$ as a function of $x$ and satisfies
$$\forall \alpha,\beta\in\N^d,\;\; \exists C_{\alpha,\beta}>0,\;\;\forall (x,\xi)\in\Omega\times\R^d,\;\;
\left| \partial_x^\alpha\partial_\xi^\beta a(x,\xi)\right|\leq C_{\alpha,\beta} (1+|\xi|)^{m-|\beta|}.$$
Then, the operator $a(x,D)$, defined in the Weyl quantization by
$$\forall f\in L^2(\R^d),\;\;a(x,D)f(x)=(2\pi)^{-d}\int a\left({x+y\over 2},\xi\right) f(y)\,{\rm e}^{i\xi\cdot (x-y)} dy\,d\xi$$
maps $H^s(\Omega)$ into $H^{s-m}(\Omega)$ (see \cite{AG}).
We will also assume that for $a\in{\mathcal A}_i^m$, there exists a function $a_\infty(x,\xi)$ homogeneous of degree $0$ such that for
\begin{equation}\label{def:ainfini}
\forall \omega\in{\bf S}^{d-1},\;\;\lim_R,{\infty}R^{-m}a(x,R\omega)= a_\infty(x,\omega).
\end{equation}
The symbols of ${\mathcal A}_i^m$  are
called  interior symbols because they are compactly supported inside $\Omega$. Such symbols are of no help for studying the behaviour of $(f_n)_{n\in\N}$ close to $\partial \Omega$: one then use  tangential symbols which are defined in the Appendix.
It is easy to convince oneself that the limits of quadratic quantities $\left(a(x,D) f_n,f_n\right)$ (for $(f_n)_{n\in\N}$ u.b. in $L^2(\Omega)$ and $a\in{\mathcal A}_i^0$) only depends on the function $a_\infty$. Then, following~\cite{G1} and~\cite{T}, it is possible to prove that these limits are characterized by a positive Radon matrix-valued measure $\mu$ on  $\Omega\times {\bf S}^{d-1}$ such that, up to the extraction of a subsequence,
 $$\left(a(x,D) f_n\;,\;f_n\right)\td_n,{\infty} \langle a_\infty ,\mu\rangle.$$
 Such a  measure $\mu$ is called a microlocal defect measure of the family $(f_n)_{n\in\N}$.

$ $

Let us come back to the  sequences $(\partial_t u_n)_{n\in\N}$, $(\nabla_x u_n)_{n\in\N}$ and $(\theta_n)_{n\in\N}$ that we need to study simultaneously. In that purpose, we set
\begin{equation} \label{def:Un}
U_n=\begin{pmatrix}\partial_tu_n\\   \nabla_x u_n\cr \theta_n\\\end{pmatrix}
\end{equation}
and the energy density $e_n$ writes
$$e_n(t,x)=|U_n(t,x)|^2_{\C^{d+2}}.$$
We shall
 consider quadratic quantities $\left(a(x,D) U_n(t),U_n(t)\right)$ where the symbol $$a(x,\xi)=\left(a_{i,j}(x,D)\right)_{i,j}$$
 is a  $(d+2)\times (d+2)$ matrix of interior symbols of ${\mathcal A}_i^m$. Then,  a microlocal defect measure  of $(U_n(t))_{n\in\N}$ will be matrix-valued and will describe  up to the extraction of a subsequence, the limits of the quantities
 $\left(a(x,D) U_n(t)\;,\;U_n(t)\right)$. Of course,
 the $t$-dependence of these measures is an issue by itself. Therefore, since   $U_n\in L^2_{loc}(\R,L^2(\Omega))$, we  test $(a(x,D)U_n(t),U_n(t))$ against smooth compactly supported functions of the variable $t$ and consider the limits of $$I(a,\theta):=\int\theta(t)(a(x,D)U_n(t),U_n(t))dt.$$  A microlocal defect measure $M$ of $(U_n(t))_{n\in\N}$ is a positive Radon measure on the set~$\R^+\times\Omega\times{\bf S}^{d-1}$ such that, up to a subsequence,  we have: $\forall \theta\in{\mathcal C}_0^\infty(\R)$
 $$I(a,\theta)\td_n,{+\infty} \int_{\R\times\Omega\times{\bf S}^{d-1}} \theta(t) {\rm tr} \left(a_\infty(x,\omega) M(dt,dx,d\omega)\right).$$
 In particular, $\forall \theta\in{\mathcal C}_0^\infty(R)
 ,\;\;\forall \chi\in{\mathcal C}_0^\infty (\Omega)$,
 $$ \int\theta(t) \chi(x) e_n(t,x)dxdt\td_n,{\infty} \int_{\R\times\Omega\times{\bf S}^{d-1}}\theta(t) \chi(x) {\rm tr}\, M(dt,dx,d\omega).$$
The  matrix-valued measure $M(t,x,\omega)$ is positive in the sense that its diagonal componants $m_{i,i}$ are positive Radon measure and its off-diagonal componants $m_{i,j}$ are absolutely continuous with respect
to $m_{i,i}$ and $m_{i,j}$. Using the special form of the componants of the vector $U_n(t)$, one can write
 $$M(t,x,\omega)=\begin{pmatrix}
 m_1(t,x, \omega ) & 	m_{0,1}(t,x, \omega)\omega  & 	m_{0,2}(t,x, \omega)\\
  \overline{ m_{0,1}}(t,x,\omega)^t\omega &  m_0(t,x,\omega)\omega\otimes\omega &  m_{1,2}(t,x,\omega)\omega \\
\overline{  m_{0,2}}(t,x, \omega) & \overline{	m_{1,2}}(t,x, \omega )^t\omega  & 	\nu_0(t,x, \omega )
  \end{pmatrix}$$
where
\begin{itemize}
\item $m_1(t,x,\omega)$, $\nu_0(t,x,\omega)$ are the microlocal defect measures of $(\partial_t u_n)_{n\in\N}$ and of $(\theta_{n})_{n\in\N}$ respectively,
\item  $m_0(t,x,\omega)\omega_i\omega_j$ is the joint measure of $(\partial_{x_i}u_{n})_{n\in\N}$ and $(\partial_{x_j} u_{n})_{n\in\N}$,
\item  $m_{0,1}(t,x,\omega)\omega_j$ is the joint measure of $(\partial_tu_{n})_{n\in\N}$ and $(\partial_{x_j}u_{n})_{n\in\N}$,
\item  $m_{0,2}(t,x,\omega)$ is the joint measure of $(\theta_{n})_{n\in\N}$ and $(\partial_t u_{n})_{n\in\N}$,
\item  $m_{1,2}(t,x,\omega)\omega_j$ is the joint measure of $(\theta_{n})_{n\in\N}$ and $(\partial_{x_j}u_{n})_{n\in\N}$.
\end{itemize}
By ``joint measure'' of two u.b. families of $L^2(\Omega)$,  $ (f_n)_{n\in\N}$ and $(g_n)_{n\in\N}$, we mean a measure which  describes the limit up to extraction of a subsequence of quantities $(a(x,D)f_n,g_n)_{n\in\N}$ for symbols $a$ of ordre $0$. In the following, we assume that $(U_n(0))_{n\in\N}$ has only one microlocal defect measure $M(0,x,\omega)$.

\subsection{Analysis of the thermoelasticity operator}
\ni Let us now come back to our system \aref{eq:thermo}  that we rewrite as
$$i\partial_t U_n=P(x,D) U_n$$
 where $U_n$ is defined in~\aref{def:Un}
and
$$P(x,D)=\begin{pmatrix}0 & i^t\nabla & -i\nabla\cdot(\gamma\cdot)\\
i\nabla & 0 & 0 \cr -i\gamma\cdot\nabla & 0 &
i\nabla\cdot (B(x)\nabla\cdot)\\\end{pmatrix}.$$
We first study  the eigenspaces of the matrix $P(x,\xi)$ which is not self adjoint if $B\not=0$. However, we have the following proposition:

\begin{prop}\label{prop:edecomp} There exists $R_0>0$ such that for $|\xi|\geq R_0$, the matrix $P(x,\xi)$ diagonalize with smooth eigenvalues and smooth eigenprojectors. Besides, the kernel
 of $P(x,\xi)$ is of dimension greater or equal to~$d-1$. \\
 Denote by    $\lambda_-$, $\lambda_0$ and $\lambda_+$ the three eigenvalues of $P(x,\xi)$ which are not identically~$0$,  and by  $\Pi^-$, $\Pi^0$ and $\Pi^+$ the smooth associated eigenprojectors. Then, for all $R>R_0$ and for $\chi\in{\mathcal C}^\infty(\R^d)$ such that $\chi(\xi)=0$ for $|\xi|\leq {1/2}$, $\chi(\xi)=1$ for $|\xi|\geq 1$ and $0\leq \chi\leq 1$, we have  in ${\mathcal D}'(\Omega)$,
 \begin{equation}\label{eq:edecomp}
 |U_n(t,x)|^2_{\C^{d+2}}= \sum_{k\in\{0,-,+\}} \left|                                                                                                                                                                                                                                                                                           \Pi^k(x,D)\chi\left({D\over R}\right)U_n(t,x)\right|^2_{\C^{d+2}}+o(1).
 \end{equation}
\end{prop}

This proposition is a consequence of Propositions~\ref{prop:symbol} and~\ref{prop:eigenprojector} below where we study the asymptotics of $P(x,\xi)$ for large $\xi$ (which differ whether $\left(x,{\xi\over|\xi|}\right)\in \Lambda$ or not).
 Each eigenprojector $\Pi^k(x,\xi)$  ($k=0,-,+$) of $P(x,\xi)$ characterizes a mode and for each mode, we will analyze the microlocal defect measure of the  componant $(\Pi^k(x,D)U_n)_{n\in\N}$.

\subsection{Propagation and damping for weakly degenerated $(B,\gamma)$}
Let us precise first our assumptions.

\begin{assumption}\label{assumption}We say that the pair $(B,\gamma)$ is weakly degenerated if $B$ and $\gamma$ satisfy the following conditions.
\begin{enumerate}
\item There exists a hypersurface $\Sigma$ of $\R^d$ such that  $\{{\rm det}\, B(x)=0\}\subset \Sigma$ and  for all $(x,\omega)\in\Lambda$, then $\omega$ is transverse to $\Sigma$ in $x$.
\item  There exists $\tilde\gamma\in{\mathcal C}^\infty(\Omega, \R^d)$ such that $\gamma(x)=B(x) \tilde\gamma(x).$
\end{enumerate}
\end{assumption}

\begin{example} Suppose $d=2$, $\Omega=\R$, $  \gamma(x)=e_2$, $B(x)=\begin{pmatrix}	b(x_1) & 0 \\ 0 & 	1 \end{pmatrix}$. Suppose
that $b(y)=0$ if and only if $y=0$. Then, $$\Lambda = \{(0,y),(\pm1,0)),	y \in \R\}$$ and one can check that $(B,\gamma)$ is weakly degenerated.
\end{example}

The sequences
$$v_{n,0}^\pm={1\over\sqrt 2} \left( u_{n,1}\pm { D\over |D|} \cdot \nabla_x u_{n,0}\right)$$
have only one  microlocal defect measure $\mu_0^\pm$ given by
 $$\mu_0^\pm(x,\omega)=m_1(0,x,\omega)+m_0(0,x,\omega)\pm {\rm Re}\left(m_{0,1}(0, x,\omega)\right).$$

\begin{theorem}\label{theo:nondeg} Suppose that $(B,\gamma)$ is weakly degenerated, that
 ${\rm Supp}(\mu_0^\pm)\subset \Lambda^c$ and that there exists
 $\tau_0\in\R^{*+}$  such that for $(x,\omega)\in {\rm Supp}(\mu_0^\pm)$  and for  $t\in[0,\tau_0]$,
$(x\pm t\omega ,\omega)\in{\Omega}$.  Then, for
 $\chi\in{\mathcal C}_0^\infty([0,\tau_0])$ and $\phi\in{\mathcal C}_0^\infty (\{{\rm det} B(x)\not=0\})$,
\begin{eqnarray*}
\lim_{n},{+\infty} \int_{\R^+\times \Omega} \chi(t)\phi(x) e_n(t,x) dt dx  & = &  \int_{\R^+\times\Omega\times{\bf S}^{d-1}}
\chi(t)\phi(x)d\mu^+_t(x,\omega)dt\\
 &+ & \int_{\R^+\times\Omega\times{\bf S}^{d-1}}
\chi(t)\phi(x)d\mu^-_t(x,\omega)dt
\end{eqnarray*}
where for all $a\in{\mathcal C}_0^\infty(\Omega\times {\bf S}^{d-1})$,
\begin{equation}\label{eq:formule}
\langle a,\mu_t^\pm\rangle= \int_{\Omega\times{\bf S}^{d-1}} a\left(x\pm t\omega,\omega\right){\rm Exp} \left[-\int_0^{t}{\left(\gamma(x\pm\sigma\omega)\cdot \omega\right)^2\over B(x\pm\sigma\omega)\omega\cdot\omega}d\sigma\right]d\mu_0^\pm(x,\omega).
\end{equation}
\end{theorem}

\begin{rem} 1) Note that even though the support of $\phi$ does not intersect the set $\{{\rm det} B(x)=0\}$, the trajectories $x\pm t\omega$ which reach the support of $\phi $ for $t\in[0,\tau_0]$ may pass through it.  \\
2) Besides,
because of  (1) in Assumption~\ref{assumption}, a trajectory $x+ s\omega$ crosses~$\Sigma$ at a finite number of times $0< t_1< \cdots<t_N\leq \tau_0$. Assuming $B(x+t_j\omega)\omega=0$,  by (2) of Assumption~\ref{assumption}, for $j\in \{1,\cdots,N\}$, there exists $c_j\in\R$ such that
$$ {\left(\gamma(x + s\omega)\cdot\omega\right)^2\over \omega\cdot B(x +s\omega)\omega}\sim  c_j (s-t_j)\;\; as \;\; s\sim t_j.$$
This implies that the integral in~\aref{eq:formule} is well defined.
\end{rem}

This result is proved in Section~\ref{subsec:theonondeg}.
The measures $\mu_t^\pm$ contain the part of the energy corresponding to the projection of $U_n$ on the $\pm$-mode. There is no contribution of the $0$-mode (which corresponds to the temperature) because of the smoothing effect of the Heat equation.
Note finally  that the damping in~\aref{eq:formule} can be $0$ if for all times $t\in[0,\tau_0]$, we have $\gamma(x\mp t\omega)\cdot\omega=0$.
We refer to section~\ref{subsec:theonondeg} for a discussion of what happens when (2) fails in Assumption~\ref{assumption}: all the energy may be damped in finite time (see Remark~\ref{rem:totaldamp}).

$ $

Finally, in the Appendix, we shortly discuss what happens close to the boundary
 under the following assumptions:
\begin{assumption}\label{ass2}
\begin{enumerate}
\item The rays of geometric optics have no contact of infinite order with the tangent to $\partial\Omega$.
\item There exists a compact $K$ such that $\{{\rm det}\, B(x)=0\}\subset K\subset \Omega$.
\item For all $x\in\partial\Omega$, $\gamma(x)\in T_x(\partial\Omega)$
\end{enumerate}
\end{assumption}
If  Assumption~\ref{ass2} holds, then one can use the methods developed in~\cite{GL} and the papers~\cite{B}, \cite{BL} or~\cite{L1} for the analysis at the boundary of microlocal measures of a family of solutions of a wave equation. One can prove the propagation of the energy along the generalized bicharacteristics as defined in~\cite{MS}. We shortly explain in Appendix~\ref{sec:boundary} how to reduce to the analysis of~\cite{BL}, without writing in details the arguments which are exactly the same ones.

\subsection{ Distorted propagation in an open set included in $\{B(x)=0\}$} We
suppose now that $B(x)=0$ in $\widetilde\Omega$ an open subset of $\Omega$. Then,
 the symbol $P(x,\xi)$ is self-adjoint  on $\widetilde\Omega\times \R^d$ and the method of~\cite{GMMP}  can be adapted with straightforward modifications.  We set
$$c(x,\xi)=\sqrt{(\gamma(x)\cdot\xi)^2+|\xi|^2},$$ and we consider
$\tilde \mu ^\pm_0$ and $\tilde\nu_0$ the microlocal defect measures of the sequences
$$\displaylines{\widetilde v_n^\pm={1\over\sqrt 2}\left(\pm u_{n,1} - iW(x,D) |D| u_{n,0}+N(x,D)\theta_{n,0}\right)\cr
\widetilde\theta_n= iN(x,D)|D|u_{n,0} +W(x,D)\theta_{n,0}\cr}$$
where $W(x,D)$ and $N(x,D)$ are pseudodifferential operators of order $0$ of symbols
$$W(x,\xi)= {|\xi|\over  c(x,\xi)} \, \chi(\xi/R_0)\;\;{\rm and}\;\; N(x,\xi) ={\gamma(x)\cdot\xi  \over c(x,\xi)}\, \chi(\xi/R_0);$$
 the function $\chi$ is smooth and satisfies $\chi(\xi)=0$ for $|\xi|\leq 2 $ and $\chi(\xi)=1$ for $|\xi|\geq 4$ and $R_0$ is chosen as  in Proposition~\ref{prop:edecomp}. Note that the measures  $\tilde \mu ^\pm_0$ and $\tilde\nu_0$ do not depend on the cut-off function $\chi$, besides, they satisfy
\begin{eqnarray*}
\tilde\nu_0   & = &   {1\over c(x,\omega)^2} \left( M(x,\omega)
\begin{pmatrix}0\\ \left(\gamma(x)\cdot\omega\right)\,\omega \\ 1 \end{pmatrix}
\,\Biggl|\,
\begin{pmatrix}0\\\left( \gamma(x)\cdot\omega\right)\,\omega \\ 1 \end{pmatrix}\right)_{\C^{d+2}},
\\
 \tilde m^\pm_0 & = & {1\over 2 c(x,\omega)^2} \left( M(x,\omega)
\begin{pmatrix}\pm c(x,\omega) \\- \omega \\ \gamma(x)\cdot\omega\, \end{pmatrix}
\,\Biggl|\,
\begin{pmatrix}\pm c(x,\omega) \\- \omega \\ \gamma(x)\cdot\omega\, \end{pmatrix}\right)_{\C^{d+2}}.
\end{eqnarray*}
Before stating the result, let us introduce a notation. Define
$$H_c(x,\xi)= \nabla_\xi c(x,\xi)\cdot \nabla_x -\nabla_x c(x,\xi) \cdot \nabla_\xi,$$
the Hamiltonian vector field associated with $c$ and $H_c^\infty$ the vector field  induced by $H_c$ on $S^*\Omega$.

\begin{theorem}\label{theo:deg}
Suppose that  ${\rm Supp}(\widetilde \nu_0)\subset\widetilde\Omega$, ${\rm Supp}(\widetilde \mu_0^\pm)\subset \widetilde \Omega$ and that there exists
 $\tau_0\in\R^{*+}$  such that  for all $x\in\widetilde\Omega$ and $\omega\in{\bf S}^{d-1}$, the projection on $\widetilde \Omega$ of the  integral curve of $H^\infty_c$ issued from $(x,\omega)$ stays in $\widetilde \Omega$  on the time interval~$[0,\tau_0]$.
 Then, for
 $\chi\in{\mathcal C}_0^\infty(\R)$ and $\phi\in{\mathcal C}_0^\infty (\widetilde\Omega)$,
\begin{eqnarray*}
\lim_{n},{+\infty} \int_{\R^+\times \Omega} \chi(t)\phi(x) e_n(t,x) dt dx   &= & \int_{\R^+\times\Omega\times{\bf S}^{d-1}}
\chi(t)\phi(x)d\widetilde \mu^+_t(x,\omega)dt\\
&&+  \int_{\R^+\times\Omega\times{\bf S}^{d-1}}
\chi(t)\phi(x)d\widetilde \mu^-_t(x,\omega)dt \\
&&+  \int_{\R^+\times\Omega\times{\bf S}^{d-1}}
\chi(t)\phi(x)\d\widetilde\nu_0(x,\omega)dt
\end{eqnarray*}
where $\widetilde \mu_t^\pm$ satisfy the transport equations
$\partial_t \widetilde \mu_t^\pm \mp H_c^\infty(x,\omega) \widetilde \mu_t^\pm=0 $
with initial data~$\widetilde \mu_0^\pm$.
\end{theorem}

\begin{rem}
Let us call distorted bicharacteristic curves the trajectories of $S^*\Omega$  associated with~$H_c^\infty$:
the transport equation for $\widetilde \mu^\pm$ implies that the energy propagates along these trajectories.
\end{rem}

Note that inside $\widetilde\Omega$, one cannot separate the contribution to the energy density of $(u_n)$ and  of $(\theta_n)$. The measures $\widetilde \mu^\pm_t$ depend of the value at time $t=0$ of  both quantities $(u_n)$ and $(\theta_n)$.

 These results call for further works: it would be interesting to know what happens at the boundary of $\widetilde\Omega$ and how transitions occur between the two regimes.  It would be also interesting to know whether this result in $\{B(x)=0\}$ extend to the set $\Lambda$ provided $\Lambda$ is invariant by the distorted rays. The following example show such a situation.

\begin{example} Suppose $d=2$, $\Omega=\R$, $  \gamma(x)=e_1$, $B(x)=\begin{pmatrix}	b(x_2) & 0 \\ 0 & 	1 \end{pmatrix}$. Suppose
that $b(y)=0$ if and only if $y=0$. Then, $$\Lambda = \{(x,\omega),\;\exists  y\in\R, \; x=(y,0),\; \omega=(\pm1,0)\}$$
is invariant by the distorted bicharacteristic curves issued from from points of $\Lambda$:
$$\left(x^\pm_s ,\xi^\pm_s\right)	=((y\pm s\sqrt 2,0),(\pm1,0)),\;\; s\in\R.$$
Note that, in that case, the distorted trajectories issued from points of $\Lambda$ coincide with the usual ones; however, they are described with different speed. Note also that~(2) of Assumptions~\ref{assumption} is not satisfied here.
\end{example}

\subsection{Organization of the paper} The main part of the article consists in the analysis of the m.d.m.s associated with the families $\left(\Pi^k(x,D) U_n \right)$  (for $k\in\{0,+,-\}$) where the functions $\Pi^k$ are defined in Proposition~\ref{prop:edecomp}.
We begin in
Section~\ref{sec:symbol} by studying the symbol $P(x,\xi)$, which allows to prove Proposition~\ref{prop:edecomp}. Then, in  Section~\ref{sec:proof},  we prove Theorems~\ref{theo:nondeg} and~\ref{theo:deg}; they rely on the analysis of the propagation of the m.d.m.s associated with the sequences $\left(\Pi^\pm(x,D) U_n\right)$ which is the object of Section~\ref{sec:propagation}.
Finally,  the Appendix~\ref{sec:boundary} is devoted to a discussion  of  the reflexion of the measures on the boundary.

%%%%%%%%%%%%%%%%%%%%%%%%%%%%%%%%%%%%%%%%%%%%%%%%%

\section{Analysis of the symbol of $P(x,D)$}\label{sec:symbol}

In this section we  analyze the properties of the matrix $P(x,\xi)$. The main interest of Weyl quantization is that the symbol of a self-adjoint operator is real-valued. We  denote  by $ \sigma(A)$ the symbol of an operator $A$ and we have in particular
\begin{eqnarray}
 \nonumber
 \sigma(\nabla)  & = &   i\xi, \\
 \label{def1}
 \sigma(\gamma(x)\cdot \nabla) & = & i\gamma(x)\cdot\xi-{1\over 2} \nabla\cdot\gamma(x),\\
 \label{def2}
 \sigma(\nabla\cdot\left(B(x)\nabla\right))  &  = &   -B(x)\xi\cdot\xi+b_0(x),
\end{eqnarray}
where
\begin{equation}\label{def:b0}
\displaystyle{b_0(x)= -{1\over 4} \sum_{1\leq j,k\leq d}\partial_{x_j,x_k}^2 B_{jk}(x)}.
\end{equation}

 Observe that, if $d=1$,  the function $b_0$ has a sign on $\Lambda$. Indeed, if $d=1$, the points $(x,\omega)\in\Lambda$ correspond to values $x$ which are minima of $B(x)$ and in this case $b_0\leq 0$. However, in higher dimension, the function $b_0(x)$ can be positive or negative indifferently as shows the following example: choose $d=2$ and $B(x)$ such that $\Lambda=\{((0,y),(0,\pm 1),\;\;y\in\R\}$ with
$$B(x)=\begin{pmatrix} x_1^4 & 0 \\ 0 & 1 + x_2^3 \end{pmatrix}+O(|x|^5)$$
close to $(0,0)$. Then we have $$b_0(x)=-{1\over 4}\left(12x_1^2 + 6x_2 \right) +O(|x|^3).$$
Therefore, $b_0(0,y)=-{3\over 2} y+O(|y|^3)$ and $b_0$ changes of sign on $\Lambda$.
Note also that if $B(x)=0$ in an open subset $\widetilde\Omega$ of $\Omega$, then $b_0=0$ in $
\widetilde\Omega$.\\
 The eigenvalues of the matrix $P(x,\xi)$ satisfies the following proposition:

\begin{prop}\label{prop:symbol} The matrix $P(x,\xi)$ has a kernel of dimension greater or equal to $d-1$. Moreover, $P(x,\xi)$ have three smooth eigenvalues  $\lambda_-$, $\lambda_0$ and $\lambda_+$ with :\\
{\bf 1)} If $\bigl(x,{\xi\over\mid\xi\mid}\bigr)\notin\Lambda$, then ${\rm dim}\left[{\rm Ker}\,P(x,\xi)\right]=d-1$ and
\begin{eqnarray}\label{eq:lambda0}
\lambda_0(x,\xi) & = &  - i\xi\cdot B(x)\xi +i b_0(x) +i{(\gamma(x)\cdot\xi)^2\over \xi\cdot B(x)\xi}  +O(|\xi|^{-1}),\\
\nonumber
\lambda_\pm(x,\xi) & = & \pm \beta(x,\xi)+ i\alpha(x,\xi),\end{eqnarray}
with \begin{equation}
\label{eq:alpha} \beta(x,\xi)=|\xi|+O(1),
\;\; \alpha(x,\xi)=-{1\over 2}{ (\gamma(x)\cdot\xi)^2\over
\xi\cdot B(x)\xi}+O(|\xi|^{-1}).
\end{equation}

\noindent{\bf 2)} If $\bigl(x,{\xi\over\mid\xi\mid}\bigr)\in\Lambda$,  then, if $b_0(x)=0$,  ${\rm dim}\left[{\rm Ker}\,P(x,\xi)\right]=d$ and if $b_0(x)\not=0$,  ${\rm dim}\left[{\rm Ker}\,P(x,\xi)\right]=d-1$. Moreover, 
\begin{eqnarray}\label{eq:surlambda}
\lambda_0(x,\xi)  &=&  ib_0(x){|\xi|^2\over c(x,\xi)^2}+O(|\xi|^{-1}),\\
\nonumber
\lambda_\pm(x,\xi)  &=&\pm \beta(x,\xi)+ i\alpha(x,\xi)
\end{eqnarray}
with
\begin{equation}\label{eq:alpha'}
\beta(x,\xi) =  c(x,\xi)+O(1) ,\;\; \alpha(x,\xi) = { 1\over 2}b_0(x)
{(\gamma(x)\cdot\xi)^2\over c(x,\xi)^2}+O(|\xi|^{-1}),
\end{equation}
and $
c(x,\xi)=\sqrt{(\gamma(x)\cdot\xi)^2 +|\xi|^2}$.
\end{prop}

The modes $\pm$ (corresponding to the eigenvalues $\lambda_\pm$)  give the wave feature of the equation. The speed of propagation is characterized by the function $\beta$ and the function $\alpha$ corresponds to the damping.  Note that in the second case, the speed of propagation are  distorted by comparison with the initial wave operator $\partial_t^2-\Delta$. Outside $\Lambda$, the eigenvalue $\lambda_0$ encounters of the Heat aspect. \\

\ni{\it Proof}:  We write $ P(x,\xi)=iQ(x,\xi)$ and for simplicity we work with $Q(x,\xi)$.  For $p,q\in\N^*$, we
denote by  $0_{p,q}$ the  $p\times q$ matrix with all coefficients equal to $0$. We have
$$Q(x,\xi)=\begin{pmatrix}
0  & i\xi&- \overline{k(x,\xi)} \\
i\,^t\xi & 0_{d,d} & 0_{d,1}\cr k(x,\xi)  & 0_{1,d}  & b(x,\xi)  \end{pmatrix},$$ where, in view of~\aref{def1} and~\aref{def2},  $b=b_2+b_0$ with $b_2=-B(x)\xi\cdot\xi $ real-valued
 and $k=ik_1+k_0$ with $$k_1=-\gamma(x)\cdot\xi ,\;\;k_0= {1\over 2} \nabla\cdot\gamma(x).$$
 The vector $(x,Y,y)\in \R\times \R^d\times \R$ is an eigenvector of $Q$ for the eigenvalue $\nu $ if and only if
 $$\left\{\begin{array}{lllllll} (1)\qquad \nu x & - & i \xi \cdot Y & +& \overline k y& = & 0,\\
(2)\qquad ix\xi & - & \nu Y& & & = & 0,\\
(3)\qquad k x&&&+& (b-\nu)y&=&0.\end{array}\right.$$

Let us suppose first that $\nu=0$, then for  $\xi\not=0$, equation $(2)$ gives $x=0$. For $b\not=0$, equation $(3)$ then implies $y=0$ and $(1)$ gives $\xi\cdot Y=0$. Therefore, for large $\xi$ the kernel of $Q(x,\xi)$ is of dimension at least $d-1$ and contains the smooth vectors
$$V_j(\xi)=\left(0,e_j(\xi),0\right)$$
where $(e_j(\xi))_{1\leq j\leq d}$ is a smooth basis of the hyperplane orthogonal to $\xi$ (such a basis exists for $\xi\not=0$).

$ $

Let us now suppose $\nu\not=0$. Equation $(2)$ implies that $Y$ is colinear to $\xi$ and
$x= -i\nu r$ where $Y=r\xi$. Equations $(1)$ and $(3)$ become
$$\left\{\begin{array}{ccccc}
\overline k y & - & i(\nu^2+|\xi|^2)r & = & 0,\\
(b-\nu)y & - & i\nu k r & = & 0 .
\end{array}\right.$$
Therefore, the non zero eigenvalues are the roots  of  the real-valued polynomial
$$f(X)=-X^3+X^2 b-X(|\xi|^2 +|k|^2)+b|\xi|^2.$$
It is easy to see that for large $\xi$, $f$ has one real-valued root $\nu_0$ (with $\nu_0\not=0$ for $b_0\not=0$) and two conjugated complex-valued roots $\nu_+$ and $\nu_-$. These three roots of the polynomial $f$  are simple and thus smooth; they give three smooth eigenvalues of $Q$. Consider three associated eigenvectors $V_0$, $V_+$ and $V_-$, they  are independent from the vectors $V_j(\xi)$ defined above. Therefore, we are left with a basis of eigenvectors: the matrix $Q(x,\xi)$ diagonalize.

$ $

Let us now study more precisely the asymptotics of the eigenvalues.

\noindent {\bf 1)}
 Suppose $(x,\xi)\notin\Lambda$, $\xi$ large and denote by $X_1<X_2<0$ the two negative roots  of  $f'(X)$. We have
 $X_1={2\over 3} b_2 +O(|\xi|)$ and $X_2=O(|\xi|)$. Then, using $f'(X_2)=0$, we obtain  $f(X_2)=b(X_2^2+|\xi|^2)<0$ and deduce that $f$ has only one real-valued root $\nu_0$ with $\nu_0<X_1$. We
set $\nu_0=\phi b$ with $\phi\geq {2\over 3} +o(1)$ and
\begin{equation}\label{eq:f}
0=f(\nu_0)= b^3\phi^2(1-\phi)-\phi
b(|k|^2+|\xi|^2)+b|\xi|^2.
\end{equation}
Necessarily,
$\phi=1+r$ with
$r=-{k_1^2/b_2^2}+O\left(|\xi|^{-3}\right),$ whence $\nu_0=b -{k_1^2\over b_2}+O\left(|\xi|^{-1}\right)$ and \aref{eq:lambda0}.\\
\ni Let $\nu_\pm$ be the two other (non real-valued) roots, we set
$\nu_\pm=\alpha\mp i \beta$.
 We observe
\begin{equation}\label{eq:nu+-}
\left\{\begin{array}{l}
\nu_++\nu_-=2\alpha=b-\nu_0,\\
\nu_0(\nu_++\nu_-) +\nu_+\nu_-= 2\alpha\nu_0 +\alpha^2+\beta^2=
|k|^2+|\xi|^2,\end{array}\right.\end{equation}
whence $\alpha={k_1^2\over 2b_2}+O(|\xi|^{-1})$ and $\beta^2=|\xi|^2+O(|\xi|)$. This implies
 \aref{eq:alpha}.

$ $

\noindent {\bf 2)} Suppose now $(x,\xi)\in\Lambda$, then $b(x,\xi)=b_0(x)$. The polynomial function $f'(X)$ has no real-valued root and $f(X)$ has only one real-valued  root   $\nu_0$. Since $f(0)f(b_0)<0$, the function $\nu_0(x,\xi)=O(1)$ as $\xi$ grows and
$$\nu_0(|\xi|^2+|k|^2)=b_0|\xi|^2+O(1),$$
whence~\aref{eq:surlambda}.
Besides,  denoting as before by $\nu_\pm=\alpha\mp i\beta$ the two other roots,  \aref{eq:nu+-} gives $$\alpha={b_0\over 2} {k_1^2\over k_1^2+|\xi|^2}+O(|\xi|^{-1})\;\;{\rm and }\;\;\beta^2= k_1^2+|\xi|^2+O(1),$$ whence \aref{eq:alpha'}.
$\diamondsuit$

$ $

\noindent Let us now
describe the eigenspaces of $P(x,\xi)$ for large $\xi$.

\begin{prop}\label{prop:eigenprojector} For $\xi$ large enough, there exists a  smooth basis of $\C^{2d+2}$, $$\displaylines{(V_0,V_+,V_-,V_j), \;\;1\leq j\leq d-1\;\; such\;\; that\cr
P(x,\xi)V_k=\lambda_kV_k\;\; for \;\;\in\{0,+,-\},\cr
 P(x,\xi)V_j=0\;\; for \;\;1\leq j\leq d-1,\cr}$$ with the following
expansions\\
{\bf 1)} If $\left(x,{\xi\over|\xi|}\right)\notin\Lambda$,
\begin{eqnarray}
\label{eq:Vj} V_j&  = &  (0, e_j(\xi), 0),\;\;e_j(\xi)\cdot\xi=0, \;\;1\leq j\leq d-1,\\
\label{eq:Vpm} \; V_{\pm} & = & {i\over\sqrt 2}\left(\mp 1, \,{\xi\over|\xi|},0\right)+O(|\xi|^{- 1}),\\
\label{eq:V0} V_0 & = & (0,0,1)+O(|\xi|^{- 2}).
\end{eqnarray}
{\bf 2)} If $\left(x,{\xi\over|\xi|}\right)\in\Lambda$,
\begin{eqnarray}
\label{eq:Vj*} V_j&  = &  (0, e_j(\xi), 0),\;\;e_j(\xi)\cdot\xi=0, \;\;1\leq j\leq d-1,\\
\label{eq:Vpm*} V_{\pm} & = & {1\over\sqrt 2} \left(-1,\pm \,{\xi \over c(x,\xi)}, \mp{ \gamma(x).\xi \over c(x,\xi)}\right)+O(|\xi|^{- 1}),\\
\label{eq:V0*} V_0 & = &\left(0,{ \gamma(x).\xi \over c(x,\xi)} {\xi\over|\xi|},{|\xi|\over c(x,\xi)}\right)+O(|\xi|^{- 1}).
\end{eqnarray}
\end{prop}

\noindent Note that there exits  smooth eigenvectors but their asmptotics are discontinuous; similarly, even though they are not orthogonal projectors, their asymptotics give a decomposition of $\C^{2d+2}$ on orthogonal subspaces. Moreover, when $b_0=0$, there is an eigenvalue crossing between the eigenspace for  $\lambda_0$ which merges into the kernel of $P$. However,  there still exists a smooth basis of eigenvectors and we will take advantage of this fact in the following sections. Besides, we have the following remark.

\begin{rem}\label{rem:33}
  Assuming (2) of Assumption~\ref{assumption}, we have for $(x,\omega)\in \Lambda$, $$\gamma(x)\cdot \omega=\widetilde \gamma(x)\cdot \left(B(x)\omega\right)=0.$$
  Therefore, if $\left(x,{\xi\over |\xi|}\right)\in\Lambda$,  $c(x,\xi)=|\xi|$ and
  \begin{eqnarray*}
V_j&  = &  (0, e_j(\xi), 0),\;\;e_j(\xi)\cdot\xi=0, \;\;1\leq j\leq d-1,\\
 V_{\pm} & = & {1\over\sqrt 2} \left(-1,\pm \,{\xi \over |\xi|},0\right)+O(|\xi|^{- 1}),\\
 V_0 & = &\left(0,0,1\right)+O(|\xi|^{- 1}).
\end{eqnarray*}
\end{rem}

\noindent  Let us now prove Proposition~\ref{prop:eigenprojector}.

$ $

\noindent{\it Proof :} In the proof of Proposition~\ref{prop:symbol}, we have already built the vectors $V_j$. Besides,  we have seen that the eigenvectors of $Q(x,\xi)$ associated with $\nu_0$, $\nu_+$ and $\nu_-$ are of the form $(x,r\xi,y)$ with
$x= -i\nu r$ and
$$\left\{\begin{array}{ccccc}
(b-\nu)y & - & i\nu k r & = & 0 ,\\
\overline k y & - & i(\nu^2+|\xi|^2)r & = & 0.
\end{array}\right.$$

Let us consider first the $\pm$-modes. We have $\nu_\pm=O(|\xi|)$. Therefore, $b-\nu_\pm\not=0$ for large $\xi$ independently of the fact that $(x,\xi/|\xi|)\in\Lambda$ or not. The vectors
\begin{equation*}\label{VetK}
V_\pm(x,\xi)=\tilde r\left(\nu_\pm(\nu_\pm-b),i(\nu_\pm-b)\xi,k\nu_\pm\right),
\end{equation*}
 with $ \tilde r=\left(|\nu_\pm-b|^2(|\nu_\pm|^2+|\xi|^2)+|k|^2|\nu_\pm|^2\right)^{-1/2}$
are smooth non-zero eigenvectors associated with $\nu_\pm$.
In view of the asymptotics of Proposition~\ref{prop:symbol}, we obtain asymptotics for $V_\pm$:
\begin{enumerate}
\item In $\Lambda^c$, $\nu_\pm=\mp i|\xi| +O(1)$, $k=-i\gamma(x)\cdot\xi+O(1)$ and $b=-B(x)\xi\cdot\xi +O(1)$. Therefore
$$\displaylines{\tilde r={1\over \sqrt 2} \left(|\xi| B(x)\xi\cdot\xi\right)^{-1} \left(1+O(|\xi|^{-1})\right),\cr
V_\pm(x,\xi)=\tilde r \left[\left( \mp i|\xi|(B(x)\xi\cdot\xi),i(B(x)\xi\cdot\xi)\, \xi,\mp |\xi|\gamma(x)\cdot\xi\right)+O(|\xi|^2)\right],\cr}$$
whence~\aref{eq:Vpm}.
\item  In $\Lambda$, $\nu_\pm=\mp ic(x,\xi)+O(1)$, $b=O(1)$ and we still have $k=-i\gamma(x)\cdot\xi+O(1)$. Therefore
$$\displaylines{
\tilde r=1/(\sqrt2 \,c(x,\xi)^2) (1+O(1)),\cr
V_\pm(x,\xi)=\tilde r \left[\left(-c(x,\xi)^2, \pm c(x,\xi)\,\xi,\mp c(x,\xi)\gamma(x)\cdot\xi\right)+O(|\xi|)\right],\cr
}$$
whence~\aref{eq:Vpm*}.
\end{enumerate}

$ $

Let us consider now the $0$-mode. We have $\nu_0=O(1)$ in $\Lambda$ and $\nu_0+O(1) \in\R$ in~$\Lambda^c$. Therefore, $\nu_0^2 +|\xi|^2\not=0$ for large $\xi$ and the vector
$$\displaylines{
V_0(x,\xi) = r_0\left(-\nu_0\overline k, -i\overline k \xi, \nu_0^2+|\xi|^2\right)\cr
{\rm with}\;\; r_0=\left(|\nu_0|^2|k|^2+|\xi|^2|k|^2+\left(\nu_0^2+|\xi|^2\right)^2\right)^{-1/2}\cr}$$
is a smooth eigenvector associated with the eigenvalue $\nu_0$. We are now left with a smooth basis of eigenvectors.
Let us now study the asymptotics of this vector.
\begin{enumerate}
\item In $\Lambda^c$, $\nu_0=-B(x)\xi\cdot\xi +O(1)$, therefore
$$r_0=(B(x)\xi\cdot\xi)^{-2} \left(1+O(|\xi|^{-1})\right)\;\;{\rm and}\;\; V_0(x,\xi)=(0,0,1)+O(|\xi|^{-2}).$$
\item In $\Lambda$, $\nu_0=O(1)$, whence
$$\displaylines{r_0=1/(c(x,\xi)|\xi|)(1+O(|\xi|^{-1}),\cr V_0(x,\xi)=r_0\left[\left(0,(\gamma(x)\cdot\xi)\,\xi,|\xi|^2\right)+O(|\xi|)\right],\cr}$$
which gives~\aref{eq:V0*}.$\diamondsuit$
\end{enumerate}

 $ $

 Before concluding this section, we point out that these asymptotics imply Proposition~\ref{prop:edecomp}. Indeed, we have obtained the existence of $R_0$ such that for $|\xi|>R_0$, $P(x,\xi)$ diagonalize with smooth eigenprojectors given by the Gram matrix of the basis of eigenvectors of Proposition~\ref{prop:eigenprojector}. Asymptotically, these projectors are orthogonal and we have
 $$\Pi^k(x,\xi)=V_k(x,\xi) \otimes \overline{V_k(x,\xi)}+o(1). $$
 Since $U_n(t,x)$ goes weakly to $0$ in $L^2(\Omega)$, we have
 $$|U_n(t,x)|^2_{\C^{d+2}}=\left|\chi\left({D\over R}\right)U_n(t,x)\right|^2+o(1)\;\;{\rm in}\;\;{\mathcal D}'(\Omega).$$
 Besides for $R>R_0$, we observe that $\Pi^k(x,D) \chi(D/R)U_n(t,x)$ is well-defined  and we have in ${\mathcal D}'(\Omega)$,
 $$\left|\Pi^k(x,D)\chi\left({D\over R}\right)U_n(t,x)\right|_{\C^{d+2}}^2=
 \left| \left(V_k(x,D)\,\Bigl|\,\chi\left({D\over R}\right) U_n(t,x)\right)_{\C^{d+2}}\right|^2+o(1)$$
 We can now use the asymptotics of Proposition~\ref{prop:eigenprojector} which, combined with the weak convergence to~$0$ of $U_n(t,x)$ gives
$$\displaylines{\qquad
 \sum_{k\in\{0,-,+\}} \left| \left(V_k(x,D)|\chi\left({D\over R}\right) U_n(t,x)\right)_{\C^{d+2}}\right|^2\hfill\cr\hfill = (\partial_t u_n(t,x))^2 + \left|\nabla u_n(t,x)\right|^2 +(\theta_n(t,x))^2 +o(1) \;\;{\rm in} \;\;{\mathcal D}'(\Omega),\qquad\cr}$$
 whence Proposition~\ref{prop:edecomp}.

%%%%%%%%%%%%%%%%%%%%%%%%%%%%%%%%%%%%%%%%%%

%%%%%%%%%%%%%%%%%%%%%%%%%%%%%%%%%%%%%%%%%%%%%%%%%%%%%%%%%%%%

\section{Proof of the main results}\label{sec:proof}

The proofs of Theorem~\ref{theo:nondeg} and Theorem~\ref{theo:deg} are inspired by the method developed in~\cite{GMMP} for analyzing semi-classical measures associated with solutions of a system of p.d.e's. The proof of Theorem~\ref{theo:deg} is a direct adaptation of the results of~\cite{GMMP} in the microlocal defect measures setting while the proof  of Theorem~\ref{theo:nondeg} requires non-trivial adaptations due to the fact that $P(x,\xi)$ is not self-adjoint and  that one of its eigenvalue is a  symbol of order~$2$. Therefore, we focus on the proof of Theorem~\ref{theo:nondeg} and we let to the reader the simple adaptation of these arguments to prove Theorem~\ref{theo:deg}. Theorem~\ref{theo:nondeg} relies
 on  Propositions~\ref{prop:nu} and~\ref{prop:measure} stated in  Subsection~\ref{subsec:pre}; then the proof of Theorem~\ref{theo:nondeg} is performed in Subsection~\ref{subsec:theonondeg}.

\subsection{Preliminaries}\label{subsec:pre}  We state technical results (Propositions~\ref{prop:nu} and~\ref{prop:measure} ) that we will use in the next subsection.
 The proof of Proposition~\ref{prop:nu} is done at the end of this section while the one of Proposition~\ref{prop:measure} is postponed to Section~\ref{sec:propagation}.

 $ $

 \noindent The first result describes the evolution  of the temperature $\theta_n$:

\begin{prop}\label{prop:nu} Let $\nu$ be a microlocal defect
measure  of the sequence $(\theta_n)_{n\in\N}$ of $L^2_{loc}(\R,L^2(\R^d))$, then, as a measure on $\R^+\times\Lambda^c$,
$\nu=0$.
\end{prop}

$ $

\noindent Our second result concerns  the contribution to the energy density  of  the sequences $\left(\Pi^\pm(x,D)\chi(D/R)u_n\right)_{n\in\N}$ for $R>R_0$ and $\chi$ as in Proposition~\ref{prop:edecomp}.
We set
$$U^{\pm}_{n,R}(t,x)={\Pi}_\pm(x,D)\chi\left(\frac{D}{R}\right)U_n.$$
The sequences $(U^\pm_{n,R})_{n\in\N}$ are uniformly bounded in $L^2(\Omega)$ and their
 microlocal defect measures are matrix-valued measures independent of $R>0$. Besides, by the definition of $U_{n,R}^\pm$, these measures are of the form $\mu_\pm(t,x,\omega)\Pi^\pm(x,\omega)$ where the
  measures $\mu^\pm(t,x,\omega)$ can  be understood as the traces of $\mu_\pm(t,x,\omega)\Pi^\pm(x,\omega)$. We prove the following result.

\begin{prop}\label{prop:measure} Assume condition (2) of Assumption~\ref{assumption} is satisfied in $\Omega_1\subset \Omega$ and let $T>0$. There exists a subsequence $n_k$ and a continuous map $t\mapsto \mu_\pm(t)$ from $[0,T]$ into the set of positive Radon measures on $\Omega\times{\bf S}^{d-1}$ such that for all $t\in[0,T]$ and for all scalar symbol $a\in{\mathcal A}_i^0$, we have
\begin{equation}\label{mupm}
\left(a(x,D_x)U^\pm_{n_k,R}(t)\;|\; U^\pm_{n_k,R} (t)\right) \td_k,{+\infty} \int a_\infty(x,\omega) d\mu_\pm(t,x,\omega).
\end{equation}
Moreover,
in ${\mathcal D}'(\{t\geq 0\}\times\Omega_1\times{\bf S}^{d-1})$, we have
\begin{equation}\label{eq:mu}
\partial_t  \mu _\pm  \pm H_\beta ( \mu _\pm)-2\alpha_\infty   \mu _\pm=\nu_\pm
\end{equation}
where $\nu_\pm$ is a measure supported on $\Lambda$ absolutely continuous with respect to $\mu_\pm$ and
where for all $a\in{\mathcal A}_i^0$,
$$\int a_\infty(x,\omega)\,\d\left(H_\beta(\mu_\pm)\right) = -\int (H_\beta a)_\infty\,  d\mu_\pm =\int (\{a,\beta\})_\infty \,\d \mu_\pm.$$
\end{prop}

\noindent Proposition~\ref{prop:measure} is proved in Section~\ref{sec:propagation} below. \\
Let us now prove Proposition~\ref{prop:nu}.

$ $

\ni{\it Proof}: It is enough to prove that if $q$ is a symbol of order $0$ such that  $q(x,\omega)\in{\mathcal C}^\infty(\Lambda^c)$, then
$q(x,D) \theta_n$ goes to $0$ in $L^2_{loc}(\R,L^2(\Omega)).$
 We observe that we only need to consider large values of $\xi$. Indeed, if $\chi\in{\mathcal C}^\infty (\R^d)$, $\chi(\xi)=0$ for $|\xi|\leq 1$ and $\chi(\xi)=1$ for $|\xi|\geq 2$ with $0\leq \chi\leq 1$, we have
$$q(x,D) \theta_n=q(x,D)\chi\left({D\over R}\right)\theta_n + o(1)$$
in $L^2(\Omega)$
because $\theta_n$ goes weakly to $0$ and the operator $1-\chi\left({D\over R}\right)$ is compact. We write
$$q(x,\xi)\chi\left({\xi\over R}\right)= {1\over R} Q_R(x,\xi) \cdot \sqrt{ B(x) }\xi $$ where $Q_R$ is the vector-valued symbol of order~$-1$:
\begin{eqnarray*}
Q_R(x,\xi) & = & R\; \chi\left({\xi\over R}\right) {q(x,\xi)\over B(x)\xi\cdot \xi}\sqrt {B(x)}\xi\\
& = &\; \widetilde \chi\left({\xi\over R}\right) {|\xi| q(x,\xi) \over B(x)\xi\cdot \xi} \sqrt{B(x)}\xi
\end{eqnarray*}
with $\widetilde\chi(u)=|u|^{-1}\chi(u)$. Note that $Q_R$ is smooth since $q=0$ in a neighborhood of $\Lambda$ and because $\widetilde \chi\in{\mathcal C}_0^\infty(\R^d)$. Note also  that $Q_R$ satisfies symbols estimates uniformly with respect to $R$.  Therefore, there exists a constant $C>0$ such that for all $f\in L^2(\Omega)$,
$$ \| Q_R(x,D) f\|_{L^2(\Omega)}\leq C\;\|f\|_{L^2(\Omega)}.$$
We write
$$q(x,D)\theta_n=R^{-1} Q_R(x,D) \sqrt {B(x)} \nabla_x \theta_n+o(1)$$
in $L^2(\Omega)$ where the $o(1)$ is uniform in $R>1$ when $n$ goes to $+\infty$. As a consequence, if
 $f\in L^2(\Omega)$
 and $\psi\in{\mathcal C}_0^\infty(\R)$, we have
$$\displaylines{\qquad
\Bigl|\int \psi(t)\left(q(x,D) \theta_n(t)\;|\;f\right) dt \Bigr|\hfill\cr\hfill
 =  {1\over R} \left| \int \psi(t)\left( \sqrt {B(x)}\nabla \theta_n(t) \;|\; Q_R(x,D)^* f \right)dt \right|+o(1)\qquad\cr\hfill
  \leq    {C\over R} \,\| f\| _{L^2(\Omega)} \int \left|\psi(t)\right| \|  \sqrt {B(x)}\nabla \theta_n(t)\|_{L^2(\Omega)} dt+o(1)\qquad
\cr}$$
for some constant $C>0$ and where the $o(1)$ is uniform in $R>1$ as $n$ goes to $+\infty$.
We observe that the energy equality~\aref{eq:energy} gives that the family $\left(\sqrt{B(x)} \nabla_x \theta_n(t)\right)$ is u. b. in $L ^2\left([0,T]\times \Omega,\R^d\right)$ for all $T>0$.  Therefore, letting $n$ and $R$  go to $\infty$, we obtain the result.
$\diamondsuit$

\subsection{Proof of Theorem~\ref{theo:nondeg}}\label{subsec:theonondeg}
We suppose that $(B,\gamma)$ is weakly degenerated, that is that the  conditions (1), (2) and (3) of Assumption~\ref{assumption} are satisfied. Let $\chi$, $\phi$ be as in Theorem~\ref{theo:nondeg}. We observe that the energy density $e_n(t,x)$ writes
$$ e_n(t,x)=|\theta_n(t,x)|^2+ {1\over 2}\left|\partial_t u_n -{D\over|D|}\cdot \nabla u_n \right|^2 +  {1\over 2}\left|\partial_t u_n +{D\over|D|}\cdot \nabla u_n \right|^2.$$
Therefore, it remains to analyze the limit of each of these terms for $\phi$ supported outside $\{{\rm det} B(x)=0\}$.

$ $

\noindent First, we observe that by Proposition~\ref{prop:nu}, if $\phi\in{\mathcal C}_0^\infty\left(\{{\rm det} B(x)\not=0\}\right)$,
$$\lim_{n},{+\infty}\int_\R\int_{\Omega} \chi(t)\phi(x) |\theta_n(t,x)|^2 dtdx=\int_\R\int_{\Lambda^c} \chi(t)\phi(x)d\nu(t,x,\omega)=0.$$

\noindent Therefore, the weak limit of the energy density express only in terms of the sequences
$\partial_t u_n \pm{D\over|D|}\cdot \nabla u_n$, of which  microlocal defect measures are $\mu_\pm(t,x,\omega)$ by Remark~\ref{rem:33}.
By Proposition~\ref{prop:measure}, the measures $\mu_\pm$ satisfy
$$\partial_t \mu_\pm\pm \omega\nabla_x\mu_\pm=-{\bf 1}_{\Lambda^c}\, {(\gamma(x)\cdot\omega)^2\over B(x)\omega\cdot\omega}\mu_\pm+\nu_\pm$$
where we have used $\alpha_\infty(x,\omega){\bf 1}_\Lambda={1\over 2}b_0(x){(\gamma(x)\cdot \omega)^2\over c(x,\omega)^2}=0$  by Remark~\ref{rem:33}. \\
Let us now prove that the fact that $\omega$ is transverse to the hypersurface $\Sigma$, implies that  $\mu_\pm {\bf 1}_\Sigma=0.$ Let $f(x)=0$ be a local equation of $\Sigma$ in a subset $\Omega_2$ of $\Omega$ and let $a$ be a symbol  supported in $\Omega_2\times \R^d$ and such that $a\geq 0$ and $\xi\cdot \nabla f(x)>0$ for $(x,\xi)\in{\rm Supp} \,a$, $\xi\not= 0$. We choose a function $\theta\in{\mathcal C}_0(\R)$ such that $\theta'(0)=1$ and we use the test symbol
$$b_\delta(x,\xi)=\delta \,a(x,\xi)\theta\left({f(x)\over \delta}\right)$$ where $\delta>0$. The transport equations for $\mu_\pm$ imply
$$\int a_\infty(x,\omega)\, \omega\cdot \nabla f(x)\theta'\left({f(x)\over \delta}\right) d\mu_\pm(x,\omega)=O(\delta).$$
By letting $\delta$ go to $0$, we obtain $\int_{f(x)=0} a_\infty(x,\omega) \omega\cdot \nabla f(x) d\mu_\pm(x,\omega)=0$, whence $\mu_\pm {\bf 1}_\Sigma=0$ on the support of $a$. By this way (inspired from~\cite{DFJ}), we finally obtain $\mu_\pm{\bf 1}_\Lambda=0$ since $\Lambda\subset \Sigma\times \R^d$.\\
Besides, since the measure $\nu_\pm$ is supported on $\Lambda$ and absolutely continuous with respect to $\mu_\pm$, we deduce $\nu_\pm=0$.\\
Finally, we observe that the function $F(x,\omega)=- {(\gamma(x)\cdot\omega)^2\over B(x)\omega\cdot\omega}$ extends continuously to the set~$\Omega\times {\bf S}^{d-1}$  with $F(x,\omega)=0$ on $\Lambda$ (since $\gamma(x)\cdot\omega= \tilde\gamma(x)\cdot (B(x)\omega)=O(|B(x)\omega|)$ by (2) of Assumptions~\ref{assumption}). Therefore, we can write
 $$\partial_t \mu_\pm\pm \omega\nabla_x\mu_\pm=-{(\gamma(x)\cdot\omega)^2\over B(x)\omega\cdot\omega}\mu_\pm. $$
As a conclusion, we obtain Theorem~\ref{theo:nondeg}. Indeed, take $a\in{\mathcal C}_0^\infty(\Omega\times{\bf S}^{d-1})$, $0\leq s\leq t\leq \tau_0$ and set
$$a_t(s,x,\omega)=a\left(x\pm(t-s)\omega,\omega\right){\rm Exp} \left[\int_0^{t-s}F(x\pm\sigma\omega,\omega)d\sigma\right],$$
then $a_t(t,x,\omega)=a(x,\omega)$
$${d\over ds} \langle a_t(s),\mu^\pm(s)\rangle=0$$
whence formula~\aref{eq:formule}.

$ $

Let us conclude this section by a remark.

\begin{rem} \label{rem:totaldamp}
 If (2) of Assumption~\ref{assumption} fails, then all the energy may be damped in finite time.
\end{rem}

\noindent{\it Proof}:
Suppose $\mu_0^-=0$ and $\mu_0^+=a_0\delta(x-x_0)\otimes \delta(\omega-\omega_0)$ for $(x_0,\omega_0)\notin \Lambda$ and suppose that there exists $\tau_0$ such that
$$\forall t\in[0,t_0[,\;\;  (x_0+t\omega_0,\omega_0)\notin\Lambda\;\;{\rm and}\;\;
(x_0+t_0\omega_0,\omega_0)\in\Lambda.$$
Set $y_0=x_0+t_0\omega_0$ and suppose that $dB(y_0)[\omega_0,\omega_0,\omega_0]\not=0$ and $\gamma(y_0)\cdot\omega_0\not=0$ then
$${\rm exp} \left[-\int_0^t {(\gamma(x_0+s\omega_0)\cdot\omega_0)^2\over B(x_0+s\omega_0)\omega_0\cdot\omega_0}ds\right]\td_t,{t_0} 0$$
and all the energy is damped between times $0$ and $t_0$.

%%%%%%%%%%%%%%%%%%%%%%%%%%%%%%%%%%%%%%%%%%%%%%%%%%%%%%%%%%%%%%

\section{Propagation of Microlocal defect measures}\label{sec:propagation}

In this section, we prove Proposition~\ref{prop:measure}. We  consider the $+$ mode; the $-$ mode can be  treated in the same way.  We proceed in three steps:
\begin{itemize}
\item First, we analyze the time derivative of $\left(a(x,D) U^+_{n,R}(t)\;|\;U^+_{n,R}(t)\right)$ for scalar valued symbols $a\in {\mathcal A}_i^0$ and prove that there exists a symbol $T\in{\mathcal A}_i$ such that,  uniformly in $R$ as $n$ goes to $+\infty$,
\begin{eqnarray}\label{dUdt}
\qquad{d\over dt} \left(a(x,D) U^+_{n,R}\;|\;U^+_{n,R}\right) &  = &  -\big(\{a,\beta\}(x,D)U^{+}_{n,R}|U^{+}_{n,R}\big)\\
\nonumber & & +2\big((a \alpha)(x,D)U^+_{n,R}|U^+_{n,R}\big)\\
\nonumber
& & + \left(T(x,D)\chi\left({D\over R}\right)U_n \,|\,\chi\left({D\over R}\right) U_n\right)+o(1).
\end{eqnarray}

\item We calculate precisely the  symbol $T(x,\xi)$ and show that $T\in{\mathcal A}_i^0$. Therefore, the quantity $\left({d\over dt}\left(a(x,D)U^\pm_{n,R}\;|\;U_{n,R}^\pm\right)\right)_{n\in\N}$ is uniformly bounded and by considering a dense subset of ${\mathcal A}_0^i$, Ascoli Theorem yields the existence of the continuous map $t\mapsto \mu_\pm(t)$ satisfying~\aref{mupm}.
\item Finally, we prove that for all $\psi\in{\mathcal C}_0^\infty(\R)$,
\begin{eqnarray}\label{Tinfty}
\lim_{n},{ +\infty} \int\psi(t) \left(T(x,D)\chi\Bigl({D\over R}\right)U_n(t) &|&\chi\left({D\over R}\Bigr) U_n(t)\right)dt\\
\nonumber
&  =&
\int \psi(t) a_\infty(x,\omega) d\nu_\pm(x,\omega)dt
\end{eqnarray}
where $\nu_\pm$ is a measure supported on $\Lambda$ and absolutely continuous with respect to $\mu_\pm$
\end{itemize}
At the end of these three steps, we have obtained Proposition~\ref{prop:measure}. We now detail each of these steps.

\subsection{1st. step: proof of~(\ref{dUdt})}
The family of functions $U^{+}_{n,R}$ is solution of the equation
 \begin{equation}\label{eq:U^+}
 i \partial _tU^{+}_{n,R}(t,x)=\Pi^+(x,D)\chi\left(\frac{D}{R}\right)P(x,D)U_n.
\end{equation}
Since $\nabla\chi$ is compactly supported, we have
$$\chi\left(\frac{D}{R}\right)P(x,D)=P(x,D)\chi\left(\frac{D}{R}\right)+K_R(x,D),$$
where $K_R$ is a compact operator.
Moreover we have$$\Pi^+(x,D)P(x,D)={\lambda}_+(x,D)\Pi^+(x,D)+R_1(x,D),$$
where the symbol $R_1$ will be precisely calculated in the next subsection.
So \aref{eq:U^+} becomes

\begin{eqnarray}\nonumber
i \partial_tU^{+}_{n,R}(t,x) &=& \left({\lambda}_+(x,D)\Pi^+(x,D)+R_1(x,D)\right)\chi\left(\frac{D}{R}\right)U_n\\
\nonumber & & \qquad \qquad +\Pi^+(x,D)K_R(x,D)U_n\\
\nonumber
& = & {\lambda}_+(x,D)U^{+}_{n,R}+R_1(x,D)\chi\left(\frac{D}{R}\right)U_n+\Pi^+(x,D)K_R(x,D)U_n\\
\label{eq:U^+_{n,R}}
& = &{\lambda}_+(x,D)U^{+}_{n,R}+F^{+}_{n,R}(x,D),
\end{eqnarray}
where
\begin{equation}\label{F+}
F^{+}_{n,R}(x,D)=R_1(x,D)\chi\left(\frac{D}{R}\right)U_n+\Pi^+(x,D)K_R(x,D)U_n.
\end{equation}
For  a real-valued symbol  $a\in{\mathcal A}_i^0$,  \aref{eq:U^+_{n,R}} implies
$\displaystyle{
\frac{d}{dt}(a(x,D)U^+_{n,R}|U^+_{n,R})=I_1+I_2,}$
with
\begin{eqnarray*}
I_1 & = &
\frac{1}{i}(a(x,D)\lambda _+(x,D)U^+_{n,R}|U^+_{n,R})-\frac{1}{i}(a(x,D)U^+_{n,R}|\lambda _+(x,D)U^+_{n,R}),\\
I_2 & = & \frac{1}{i}(a(x,D)F^{+}_{n,R}(x,D)|U^+_{n,R})-\frac{1}{i}(a(x,D)U^{+}_{n,R}|F^+_{n,R}).
\end{eqnarray*}
\noindent The term $I_1$ will give the transport by the vector field $H_\beta$ and the damping by~$\alpha$. The term $I_2$ is a rest term and its main contribution will be described by a symbol~$T$.

$ $

\noindent Let us study $I_1$.
$$I_1=\frac{1}{i}\big((a(x,D)\lambda _+(x,D)-(\lambda _+(x,D))^\ast a(x,D))U^+_{n,R}|U^+_{n,R}\big).$$
We recall that $\lambda _+(x,D)=\beta(x,D)+i\alpha(x,D),$ so, since we use the Weyl quantification for the symbols we have $(\lambda _+(x,D))^\ast=\beta(x,D)-i\alpha(x,D)$ and
$$\displaylines{
I_1=\frac{1}{i}\big([a(x,D),\beta(x,D)]U^{+}_{n,R}|U^{+}_{n,R}\big)+2\big((a \alpha)(x,D)U^+_{n,R}|U^+_{n,R}\big)\hfill\cr\hfill
+ \big(r_{-1}(x,D)U^+_{n,R}|U^+_{n,R}\big)\qquad\cr\hfill
=-\big(\{a,\beta\}(x,D)U^{+}_{n,R}|U^{+}_{n,R}\big)+2\big((a \alpha)(x,D)U^+_{n,R}|U^+_{n,R}\big)
+ \big(r_{-1}(x,D)U^+_{n,R}|U^+_{n,R}\big),\cr}$$
where $r_{-1}(x,\xi)$ will denote from now on a generic symbol  in ${\mathcal A}_i^{-1}.$ Using \aref{eq:alpha}, we then obtain if $n, R$ tend to infinity
\begin{equation}\label{eq:limit1}
I_1\rightarrow -\int \{a,\beta\}_{\infty}(x,\omega)d\mu_+(x,\omega)+2 \int a_\infty(x,\omega)\alpha_\infty(x,\omega)d\mu_+(x,\omega).
\end{equation}
Note that, for these terms, we do not need to integrate in time to get the convergence.

$ $

\noindent Let us now study $I_2$. We set $U_{n,R}=\chi\left(\frac{D}{R}\right)U_n$ and by~\ref{F+},
we obtain ~$I_2=I_{2,1}+I_{2,2}$ with
\begin{eqnarray}\nonumber
I_{2,1}(t)&=&\frac{1}{i}\big(a(x,D){\Pi}^+(x,D)K_RU_{n}(t)|\Pi^+(x,D)U_{n,R}(t)\big)\\
\nonumber
& & \qquad\qquad - \frac{1}{i}\big(a(x,D){\Pi}^+(x,D)U_{n,R}(t)|\Pi^+(x,D)K_RU_{n}(t)\big)\\
\nonumber
I_{2,2}(t)&= & \frac{1}{i}\big(a(x,D)R_1(x,D)U_{n,R}(t)|{\Pi}^+(x,D)U_{n,R}(t)\big)\\ & &
\nonumber
\qquad\qquad -\frac{1}{i}\big(a(x,D){\Pi}^+(x,D)U_{n,R}(t)|R_1(x,D)U_{n,R}(t)\big)
\end{eqnarray}
Since $(U_n)_{n\in\N}$  goes weakly to $0$ as $n$ goes to $+\infty$  and   $K_R$ is a compact operator, the sequence  $(K_RU_n)_{n\in\N}$ goes strongly to $0$. Therefore
\begin{equation}\label{eq:limit2}
\forall t\in \R^+,\;\;I_{2,1}(t)\rightarrow 0.
\end{equation}
Besides,
$I_{2,2}(t)=\left(T(x,D) U_{n,R}\mid U_{n,R}\right)$
with
$$T(x,D):={1\over i}\left(\Pi^+(x,D)^*a(x,D)  R_1(x,D)-R_1(x,D)^* a(x,D) \Pi^+(x,D)\right).$$

\subsection{2nd. step: the symbol of the rest term} In the following, rest terms in the symbol class ${\mathcal A}^{-k}_i$ for $k\in\N$ will be denoted by $r_{-k}$. We prove the following Lemma.

\begin{lem}\label{lem:T} Suppose  $\Omega_1$ is like in Proposition~\ref{prop:measure}. Then, in $\Omega_1$,  $T\in{\mathcal A}^0_i$ and
$$T(x,\xi)=T_1(x,\xi)+T_2(x,\xi)+r_{-1}(x,\xi)$$ with
\begin{eqnarray}\label{eq:T1}
T_1&= &- {a\over 2}
\Bigl(\Pi^+\{\Pi^+,\beta\}+\{\Pi^+,\beta\}\Pi^+ +
 \Pi^+\{\Pi^+,\beta\}\Pi^-
\\
\nonumber && +\Pi^-\{\Pi^+,\beta \}\Pi^+\Bigr)\\
 \label{eq:T2}
 T_2 & = &    {a\over 2 i}\left(\Pi^+\{\Pi^+,b\} K-K\{\Pi^+, b\}\Pi^+\right)
 \end{eqnarray}
Besides, $T_1\Pi_0$ and $\Pi_0T_1$ are symbols of ${\mathcal A}_i^{-1}$.
\end{lem}

\noindent Note that $T_2\in{\mathcal A}_i^{0}$ because $\Pi^+K, \,K\Pi^+ \in{\mathcal A}_i^{-1}$ by Remark~\ref{rem:33}.

$ $

\begin{proof}  We first calculate  $R_1$,
  recall that
$$
R_1(x,D)  =  \Pi^+(x,D)P(x,D)-{\lambda}_+(x,D)\Pi^+(x,D).$$
We write $P(x,\xi)=P_1(x,\xi)+ib(x,\xi)K,$
with
$$ P_1(x,\xi)=
  \begin{pmatrix} 0& -^t\xi& \gamma(x)\cdot\xi \\
-\xi & 0 _{d,d}& 0_{d,1} \\ \gamma(x).\xi& 0_{1,d} &
0\end{pmatrix}
\;\;{\rm and}\;\;
K=\begin{pmatrix} 0_{d+1,d+1} & 0_{d+1,1}\\
 0_{1,d+1} & 1\end{pmatrix}.$$
Since $b(x,\xi)$ is of order $2$, $R_1$ is a priori of order $1$
 and in view of $\Pi^+P=\lambda_+\Pi^+$, we have
$$R_1= {1\over 2i} \left(\{\Pi^+,P\}-\{\lambda_+,\Pi^+\}\right)+QK+r_{-1}$$
 where the  matrix-valued symbol $Q$ is in ${\mathcal A}^0_i$ and is the sum of second derivatives of $\Pi^+$ multiplied by second derivatives of $b$. More precisely and after ordering the terms of higher degrees, we write
 $$
 R_1={1\over 2} \{ \Pi^+\;,\; b\} K-{i\over 2}\left(  \{\Pi^+,\beta\}+ \{\Pi^+,P_1\}\right)+Q K+r_{-1},$$
where we have used that $\lambda_+=\beta +i\alpha$ and $\{\alpha,\Pi^+\}\in{\mathcal A}^{-1}_i$.
However, since
 $\gamma(x)\in {\rm Ran}\,B(x)$, then $\gamma(x)\cdot\xi=0$ for $\left(x,{\xi\over|\xi|}\right)\in\Lambda$ and, by Remark~\ref{rem:33}
\begin{equation}\label{case2}
K\Pi^+, \Pi^+K\in {\mathcal A}^{-1}_i.
\end{equation}
Derivations of these relations imply $\{\Pi^+,b\}K\in{\mathcal A}_i^0$ and $ QK\in{\mathcal A}^{-1}_i$. Finally, we find
 $$R_1={1\over 2}\{\Pi^+,b\} K-{i\over 2}\left( \{\Pi^+,\beta\}+ \{\Pi^+,P_1\}\right)+r_{-1}\in {\mathcal A}^0_i,$$ whence $T\in{\mathcal A}^0_i$ with
\begin{eqnarray*}
T &  =&{1\over i}\left( \Pi^+aR_1-R_1^*a\Pi^+\right)+r_{-1}\\
& = & {a\over 2i}\left(\Pi^+\{\Pi^+,b\} K-K\{\Pi^+, b\}\Pi^+\right)\\
& & -{a\over 2}
\left(\Pi^+\{\Pi^+,\beta\}+\{\Pi^+,\beta\}\Pi^+ +
 \Pi^+\{\Pi^+,P_1\}-\{P_1,\Pi^+\}\Pi^+\right)+r_{-1},
 \end{eqnarray*}
 where we have used $(\Pi^+)^*=\Pi^+ +r_{-1}$
and $$R_1^*={1\over 2} K\{\Pi^+,b\}+{i\over 2}\left( \{\Pi^+,\beta\} - \{P_1,\Pi^+\}\right)+r_{-1}.$$
We now transform the expression of the principal symbol of $T$. We write
$$P_1=P-ibK=\lambda_+\Pi^++\lambda_-\Pi^-+\lambda_0\Pi^0-ibK.$$
By Proposition~\ref{prop:symbol}, \ref{prop:eigenprojector} and Remark~\ref{rem:33}, we notice that
\begin{itemize}
\item
 in $\Lambda^c$, $\Pi^0=K+O(|\xi|^{-2})$ and $\lambda_0=ib+O(1) $,
 \item  in $\Lambda$, $\Pi^0=K+O(|\xi|^{-1})$, $\lambda_0=O(1)$ and $b=O(1)$.
 \end{itemize}
  Therefore,
$$P_1 =\lambda_+\Pi^++\lambda_-\Pi^-+r_0.$$
 In view of~\aref{eq:T2}, we deduce
 $T= T_2+\widetilde T$
 with
\begin{eqnarray*}
\widetilde T &= &-
{ a\over 2} \biggl(\left(\Pi^+\{\Pi^+,\beta\}+\{\Pi^+,\beta\}\Pi^+ \right) \\ & & \qquad+
\sum_{\ell\in\{+,-\}}\left( \Pi^+\{\Pi^+,\lambda_\ell\Pi^\ell\}-\{\overline\lambda_\ell\Pi^\ell,\Pi^+\}\Pi^+\right)\biggr)+r_{-1}\\
&= &- { a\over 2} \biggl(\left(\Pi^+\{\Pi^+,\beta\}+\{\Pi^+,\beta\}\Pi^+ \right) \\ & & \qquad+
\sum_{\ell\in\{+,-\}}\left( \Pi^+\{\Pi^+,\beta\Pi^\ell\}-\{\beta\Pi^\ell,\Pi^+\}\Pi^+\right)\biggr)+r_{-1}.
\end{eqnarray*}
Observing that $\{\Pi^\pm,\Pi^+\}, \{\Pi^+,\Pi^\pm\}\in{\mathcal A}_i^{-1}$ by Remark~\ref{rem:33} and equation~\aref{eq:T1}, we obtain
$$\displaylines{
\widetilde T  =  -{a\over 2} \biggl(\Pi^+\{\Pi^+,\beta\}+\{\Pi^+,\beta\}\Pi^+ +\sum_{\ell\in\{+,-\}}
\left(\Pi^+\{\Pi^+,\beta\}\Pi^\ell-\Pi^\ell\{\beta,\Pi^+\}\Pi^+ \right) \biggr)\cr\hfill +r_{-1},\qquad\cr
\qquad=T_1+r_{-1}\hfill\cr}$$
where we have used $\Pi^+\{\beta,\Pi^+\} \Pi^+=0$ (which comes from $(\Pi^+)^2=\Pi^+$ whence $\{\beta,\Pi^+\} =\Pi^+\{\beta,\Pi^+\} +\{\beta,\Pi^+\} \Pi^+$ and, multiplying by $\Pi^+$ on both sides, we obtain $\Pi^+\{\beta,\Pi^+\} \Pi^+=0$).
\noindent Notice that $\{\beta,\Pi^+\}\Pi^0\in{\mathcal A}^{-2}_i$ and $\Pi^0\{\beta,\Pi^+\}\in{\mathcal A}^{-2}_i$, whence $ \Pi^0T,\;T\Pi^0\in {\mathcal A}^{-1}_i$. \\
\end{proof}

\subsection{3rd. step: passing to the limit in the rest term}
 We use the following Lemma.

\begin{lem}\label{lem:A} Consider $A$ a smooth symbol of order $0$ supported  in $\Omega_1$ and such that $A=\Pi^+A \Pi^\ell$ or $A=\Pi^\ell A \Pi^+$ where $\Pi^\ell=\Pi^-$ or $\Pi^\ell =\Pi^j$, $1\leq j\leq d-1$.
  Then if $\psi\in{\mathcal C}_0^\infty(\R)$,
$$\int\psi(t) \left(A(x,D)U_n(t)\;,\;U_n(t)\right)dt\td_n,{+\infty} 0.$$
\end{lem}

Since the symbol $T_1$ is the sum of terms of the form $\Pi^+\{\Pi^+,\beta\}\Pi^\ell$ or $\Pi^\ell\{\Pi^+,\beta\}\Pi^+$ with $\ell=-$ or $\ell =j$, $1\leq j\leq d-1$, we can apply the  Lemma and we obtain
$$\displaylines{
\qquad \int \psi(t) \left(T(x,D) \chi\left({D\over R}\right) U_n(t)\,|\, \chi\left({D\over R}\right) U_n(t)\right)dt\hfill\cr\hfill
=
\int \psi(t) \left(T_2(x,D) \chi\left({D\over R}\right) U_n(t)\,|\, \chi\left({D\over R}\right) U_n(t)\right)dt+o(1)\qquad\cr}$$
as $n$ goes to $+\infty$. In view of~\ref{eq:T2}, we have
$$\displaylines{\qquad \limsup_n,{+\infty} \int \psi(t) \left(T_2(x,D) \chi\left({D\over R}\right) U_n(t)\,|\, \chi\left({D\over R}\right) U_n(t)\right)dt\hfill\cr\hfill
 \td_{R},{+\infty} \int \chi(t)  (T_2)_\infty(x,\omega) d\widetilde\nu_+(t,x,\omega)\qquad\cr}$$
where $\widetilde \nu_+$ is the joint measure of $KU_n=\theta_n$ and of $U^+_{n,R}$. The measure $\nu_+$ is absolutely continous with respect to $\mu_+$ and $\nu$. By Proposition~\ref{prop:nu}, $\widetilde \nu_+$ is supported on $\Lambda$ and we obtain~\aref{Tinfty} with
$\nu_+= (T_2)_\infty \widetilde \nu_+$ which is  absolutely continuous with respect to $\mu_+$ and supported on $\Lambda$.
We now focus on  the proof of  Lemma~\ref{lem:A}.

$ $

\begin{proof}  Let $\Pi^\ell $ be one of the projectors $\Pi^-$ or $\Pi^j$,  $1\leq j\leq d-1$ and let us consider a term of the form $\Pi^+A\Pi^\ell$; the proof is similar for the other terms. We have
$$[\Pi^+A\Pi^\ell\;,\;P]=(\lambda_\ell-\lambda_+)\Pi^+A\Pi^\ell.$$
Let us denote by $C^\ell$ the symbol of ${\mathcal A}_i^{-1}$
$$C^\ell=(\lambda_\ell-\lambda_+)^{-1}\Pi^+A\Pi^\ell.$$
We have $C^\ell\in{\mathcal A}^{-1}_i$ and, by Remark~\ref{rem:33}, we have $C^\ell K,KC^\ell\in{\mathcal A}^{-2}_i$.
 We
write
$$\displaylines{\qquad
\int\psi(t)\left(
(\Pi^+A\Pi^\ell)(x,D)U_{n,R}|U_{n,R}\right)_{L^2(\R^d)}\d t  \hfill\cr\hfill=
\int\psi(t)\left(
([C^\ell,P])(x,D)U_{n,R}|U_{n,R}\right)_{L^2(\R^d)}\d t.\qquad\cr}
$$
Besides,
$$\displaylines{\qquad
[C^\ell,P](x,D)  =  [C^\ell(x,D), -i\partial_t +P(x,D)] -{1\over 2i} \left(\{C^\ell,P\}-\{P,C^\ell\}\right) (x,D)\hfill\cr\hfill+r_{-1}(x,D).\qquad\cr}
$$
Using $P=P_1+ibK$ and $\{C^\ell,P_1\},\;\{P_1,C^\ell\}\in{\mathcal A}_i^{-1}$, we obtain,
$$\{C^\ell,P\}-\{P,C^\ell\}=i\{C^\ell,b\}K-iK\{b,C^\ell\}+r_{-1}.$$
Therefore, we have
\begin{eqnarray*}
[C^\ell,P](x,D) & = & [C^\ell(x,D), -i\partial_t +P(x,D)] -{1\over 2} \left(\{C^\ell,b\}K-K\{b,C^\ell\}\right)(x,D)\\
& & \qquad\quad\qquad\qquad\qquad\qquad\qquad\qquad+r_{-1}(x,D)\\
&=&  [C^\ell(x,D), -i\partial_t +P(x,D)]+r_{-1}(x,D)
\end{eqnarray*}
where we have used   $\{C^\ell ,b\}K=\{C^\ell K,b\}\in{\mathcal A}^{-1}_i$ and $K\{C^\ell,b\}=\{KC^\ell, b\}\in{\mathcal A}^{-1}_i$ (as a consequence of $KC^\ell,C^\ell K\in{\mathcal A}^{-2}_i$). In view of $P^*(x,D)=P(x,D)-2ib(x,D)K$, we can write
$$\displaylines{[C^\ell,P](x,D) =  C^\ell(x,D)\left(-i\partial_t +P(x,D)\right) - \left( -i\partial_t + P^*(x,D)\right) C^\ell(x,D) \hfill\cr\hfill -2ib(x,D) KC^\ell(x,D) + r_{-1}(x,D).\cr}$$
Finally, using  $(i\partial_t-P(x,D))U_{n,R}=K_RU_n$ with $K_R$ compact, we obtain
$$\displaylines{\qquad
\int\chi(t)\left(
(\Pi^+A\Pi^\ell)(x,D)U_{n,R}|U_{n,R}\right)_{L^2(\R^d)}\d t \hfill\cr\hfill
 =   -i\int\chi'(t) \left(C^\ell(x,D)U_{n,R}|U_{n,R}\right)_{L^2(\R^d)} \d t \qquad\cr\hfill
 -2i\int \chi(t) \left( b(x,D)K C^\ell(x,D)U_{n,R}|KU_{n,R}\right)_{L^2(\R^d)}\d t\\
+o(1).\qquad\cr}$$
Since $C^\ell$ is of order lower or equal to $-1$, we have
$$ \left(C^\ell(x,D)U_{n,R}|U_{n,R}\right)_{L^2(\R^d)} \td_{n,R},{+\infty} 0.$$
We use $$b(x,D)KU_{n,R} =\left (0,\cdots, 0,\chi\left({D\over R}\right)b(x,D)\theta_n\right)$$ with $b(x,D)\theta_n\in L^2_{loc}(\R,H^1(\Omega))$ by Proposition~\ref{prop:nu}. Since moreover $bKC^\ell$ is of order $0$ in $\Omega_1$, we have
 $$\limsup_n,{+\infty}\int \psi(t) \left( b(x,D)K C^\ell(x,D)U_{n,R}|KU_{n,R}\right)_{L^2(\R^d)}\d t\td_{R},{+\infty} 0,$$
whence
 $$\limsup_{n},{+\infty}\int\psi(t)\left(
(\Pi^+A\Pi^\ell)(x,D)U_{n,R}|U_{n,R}\right)_{L^2(\R^d)}\d t
\td_{R},{+\infty} 0.$$
\end{proof}

%%%%%%%%%%%%%%%%%%%%%%%%%%%%%%%%%%%%%%%%%%%%%%%%%%%%%%%%%%%%%%%
\appendix
\section{Analysis on the boundary}\label{sec:boundary}

In this appendix, our aim is to prove  that, under Assumtion~\ref{ass2}, the energy reflects on the boundary according to the laws of geometric optics. We shortly recall the arguments of \cite{BL} and explain how they apply to our setting. In all this section, we work in a neighborhood $\Omega_1$ of a point of $\partial\Omega$. We first recall in the first subsection the definition of Melrose-Sj\"ostrand compressed bundle  (see  \cite{MS} and the survey \cite{B}) and of the generalized  bicharacteristics. Then, in the second one,  we will link~\cite{BL}'s approach and ours (we also refer to \cite{B} and~\cite{D1}).  The last subsection will be devoted to the proof of the main statement of this Appendix.

\subsection{Melrose-Sj\"ostrandt compressed bundle and the generalized  bicharacteristics}

We work in space-time variables and set
$L=\R_t\times \Omega.$
We denote by~$(z,\zeta)$ the points of $T^*L$: $z=(t,x)$ and $\zeta=(\tau,\xi)$.
Then, the Melrose-Sj\"ostrandt compressed bundle to $L$ is given by
$$T^*_b L=\left(T^*L\setminus \{0\}\right) \cup \left( T^*\partial L\setminus\{0\}\right).$$
Quotienting by the action of $\R^+$ through homotheties, one obtains the normal compressed bundle to $\Omega$
$$S^*_bL=T^*_b L/ \R^+.$$
The projection
$$\pi: \left(T^*\R^{d+1}\right) _{|\overline L}\rightarrow T^*_bL$$
induces a topology on $T^*_b L$.
 On $T^*L$, we denote by $p_0$ the symbol of the wave operator  and by $\Sigma_0$ the projection on $T^*_bL$ of  its characteristic set
\begin{equation}\label{def:p0}
p_0(z,\zeta)=\tau^2-|\xi|^2,\;\;\Sigma_0=\pi\left(\{(z,\zeta),\;\;p_0(z,\zeta)=0\}\right).
\end{equation}
Locally, near a point of $\partial L$, we use normal geodesic coordinates $$(y,\eta)=(y_0,y_1,\cdots,y_d;\eta_0,\eta_1,\cdots,\eta_d)\in \R^{2d+2}$$ in an open set $V$ of $T^*\R^{d+1}$ such that $L\cap V=\{y_d>0\}\cap V$. Therefore,  the wave operator $\partial_t^2-\Delta $ writes $-\partial_{y_d}^2-R(y,D_{y'})$ where the symbol   $R(y,\eta')$ is a homogeneous polynomial of degree $2$ in $\eta'$; we denote by $r(y,\eta')$ its principal symbol.
We can now  distinguish between different sorts of points of $T^*\partial L\setminus \{0\}$: those who are not in $\Sigma_0$ and those who are in $\Sigma_0$ depending whether $R(y,\eta')<0$ or not. In the case  where $\rho\in\Sigma_0$,  there exists at most two points of $\{\tau^2=|\xi|^2\}=\{\eta_d^2=R(y,\eta')\}$ which are in $\pi^{-1} (\{\rho\})$, they correspond to the two roots of the equation $\eta_d^2=R(y_n,y',\eta')$.
 Consider $\rho\in T^*\partial L\setminus \{0\}$.
\begin{itemize}
\item If $\rho\notin \Sigma_0$, one says that $\rho$ is elliptic.
\item If ${\rm Card} \left( \pi^{-1}(\{\rho\})\cap\{\tau^2=|\xi|^2\}\right)=2$, $\rho$ is said to be hyperbolic.
\item   If ${\rm Card} \left( \pi^{-1}(\{\rho\})\cap\{\tau^2=|\xi|^2\}\right)=1$, $\rho$ is said to be glancing.
\end{itemize}
We denote by ${\mathcal H}$ (resp. ${\mathcal G}$) the hyperbolic (resp. glancing) points of $\partial L$.
We say that $\rho\in{\mathcal G}$ is
\begin{itemize}
\item non strictly glancing if $\partial_{y_d}r(\rho)\geq 0$,
\item strictly glancing if $\partial_{y_d} r(\rho)<0$,
\item diffractive if $\partial_{y_d} r(\rho)>0$,
\item glancing of order $k$ if
$$H^j_{r(y',0,\eta')} (\partial_{y_d}r_{|y_d=0}) (\rho)=0,\;0\leq j< k-2\;\;{\rm and} \; \;H^{k-2}_{r(y',0,\eta')} (\partial_{y_d}r_{|y_d=0}) (\rho)\not=0.$$
\end{itemize}
We denote by ${\mathcal G}_k$ (resp. ${\mathcal G}_d$) the set of points which are glancing of order $k$ (resp. diffractive).
The assumption that $\Omega$ has no contact of infinite order with its tangents~((1) in Assumptions~\ref{ass2}) consists in assuming
\begin{equation}\label{hyp:contact}
{\mathcal G}=\build\cup_{k\in\N^*}^{} {\mathcal  G}_k.
\end{equation}
The generalized bicharacteristic are defined by taking the Hamiltonian trajectory of~$p_0$ inside $\Omega$ and by specifying how the connection is made on the boundary.
The only problematic points are the glancing ones where the trajectory arrives tangentially to the boundary.  Indeed,
for $\rho\in{\mathcal G}$, we have $y_d=\eta_d=r(y,\eta')=0$.
Recall that the geodesic trajectories are generated by the Hamiltonian flow $H_{p_0}$
which writes in coordinates $(y,\eta)$
$$H_{p_0}=\left(-\nabla_{\eta'}r(y,\eta'),2\eta_d,\nabla_{y'}r(y',\eta'),\partial_{y_d}r(y,\eta')\right).$$
If $\rho\in{\mathcal G}_d$, $\partial_{y_d}r(y,\eta')> 0$, therefore $\eta_d$ and thus $y_d$ increase on the trajectory passing in $\rho$ and this trajectory remains in $\Omega$. On the contrary, if $\rho\notin{\mathcal G}_d$,  the coordinates on $\partial_{\eta_d}$ is nonpositive and the trajectory will leave $\Omega$: if $\eta_d$ decreases, it becomes non positive and so for $y_d$. To overcome this difficulty, one uses the vector
$$\tilde H(\rho)= H_{p_0}(\rho) -\partial_{y_d} r(y,\eta')\partial_{\eta_d}.$$
This vector has a coordinate on $\partial_{\eta_d}$ which is $0$.
One then  defines the generalized bicharacteristic as follows (see \cite{MS} or \cite{B}).

\begin{defi}
A generalized bicharacteristic is a continuous map $\gamma: \R\rightarrow T^*_bL$ such that
there exists a set $I$ of isolated points with
\begin{itemize}
\item
$\gamma (s) \in T^*L\cup{\mathcal G}\;\;for\;\;s\notin I \;\;\; and\;\;\;
\gamma(s)\in{\mathcal H}\;\;for \;\;s\in I.$
\item for $s\notin I$, $\gamma$ is differentiable with

$$\left\{\begin{array}{l}
\dot \gamma (s)=H_{p_0}\left(\gamma(s)\right)\;\;if \;\;\gamma(s)\in T^*L\cup{\mathcal G}_d,\\
\dot \gamma (s)=\tilde H(\rho)\;\;if \;\;\gamma(s)\in {\mathcal G}\setminus{\mathcal G}_d.
\end{array}\right.$$
\end{itemize}
\end{defi}

\ni It is proved in \cite{MS} that these definitions are intrinsic and that if assumption \aref{hyp:contact} is satisfied, then for $\rho_0\in T^*_bL\cap \Sigma_0$, there exists a unique generalized bicaracteristic such that $\gamma(0)=\rho_0$.

\subsection{Propagation near the boundary}
  Working in space-time variables, one first defines tangential symbols by use of the system of local normal geodesic coordinates: in an  open set ${\mathcal O}$ where we have  such coordinates $(y,\eta)$,  the function $a(y,\eta')\in{\mathcal  C}^\infty (\overline L\times \R^{d})$ is said to be a symbol of ${\mathcal A}_b^m$ if $a$ is compactly supported in  ${\mathcal O}$ in the variable $y$ and satisfies: $\forall \alpha,\;\beta \in\N^{d},
\;\exists C_{\alpha,\beta}>0$,
 \begin{equation}\label{estimationsymbolb}
 \forall (y,\eta)\in \overline L \times \R^{d+1},\;\left|\partial_y^\alpha\partial_{\eta'}^\beta \left(a(y,\eta')\right)\right|\leq C_{\alpha,\beta} (1+|\eta'|)^{m-|\beta|}
\end{equation}
This definition implies that the sets ${\mathcal A}^m_b$ depend on the choice of the open sets $({\mathcal O}_j)_{j\in J}$.
Then, one defines $${\mathcal A}^m=\{ a\in{\mathcal C}^\infty(\overline \Omega\times \R^{d}),\;\;\exists a^{i}\in {\mathcal A}^m_i,\;\; \exists a^{b}\in{\mathcal A}^m_b,\;\;a=a^b+a^i\},$$
and
 one associates with our sequence $u_n(t,x)$  its   $H^1$ space-time microlocal defect measure $\mu$ by:  up to extraction of a subsequence, for all  $q=q_b+q_i\in{\mathcal A}^2$,
$$(q(z,D_z) u_n,u_n)\td_n,{+\infty} \int_{ \overline L\times S^{d}(2)}q_\infty(z,\zeta)\d\mu(z,\zeta)$$
where $S^d(2)$ is the sphere of radius $2$ of $\R^{d+1}$ and for $(z,\zeta)\in S^d(2)$,  $$q_\infty(z,\zeta)=\lim_R,{+\infty} R^{-2} q(z,R\zeta).$$
Note
that on $S^*\Omega$, the measure $\mu$ is the usual microlocal defect measure  (or H-measure)
 as introduced in \cite{G1} and \cite{T}.
The link between this measure $\mu$ and the measures $\mu^\pm$ of Theorem~\ref{theo:nondeg} is described in the following Proposition.

\begin{prop}\label{muetmupm} In $\R\times\Omega_1\times S^d(2)$, we have
\begin{equation}\label{eq:decompmu0}
\mu(t,x,\tau,\xi)={1\over 2}\,\delta(\tau+1)\otimes \mu^+_t(x,\xi) \otimes dt+{1\over 2}\,\delta(\tau-1)\otimes \mu_t^-(x,\xi)\otimes dt.
\end{equation}
$$H_{p_0}\mu= 4 \tau \alpha_\infty\mu\;\; with \,\,\,\,\alpha_\infty(x,\omega)=-{1\over 2}{\left(\gamma(x)\omega\right)^2\over B(x)\omega\cdot\omega}. $$
\end{prop}

This proposition is proved in section~\ref{subsec:muetmupm} below.
 The next step consists in  analyzing $\mu$ on the boundary. We recall that if $\rho\in T^*\partial L$ is an hyperbolic or a glancing point, then there exists two points $\rho^\pm$  such that $\pi(\rho^\pm)=\rho$ and $\rho^\pm\in \{\tau^2=|\xi|^2\}$. These two points  are equal if $\rho\in{\mathcal G}$.  When $\rho\in{\mathcal H}$, they differ by their $\xi$-component that we will denote by $\xi^\pm$.
 Besides, as in~\cite{BL}, because of the equation satisfied by~$(u_n)_{n\in\N}$,  the sequence of the normal derivatives $(\partial_N u_n)_{n\in\N}$ is a uniformly bounded family of $L^2_{loc}(\R,L^2(\partial\Omega))$ and we denote by $\lambda$ its microlocal defect measure (we suppose that we have extracted a subsequence so that $\lambda$ is uniquely determined).
Following~\cite{BL}, one aims at calculating
$$\ell:=H_{p_0}\mu -4\tau\, \alpha_\infty\, \mu. $$
We already know that if $q$ is an interior symbol, $\langle q,\ell\rangle=0$.  For analyzing the action of $\ell$ on tangential symbols, one computes the quantity
$$I_n(t)=\left( \left[  \partial_{y_d}^2 +R(y,D_{y'}), q_b(y,D_{y'})\right] u_n(t)\;,\;u_n(t)\right)$$
for $q_b$ a tangential symbol of order $1$.

The computation of $I_n$ and the passage $n\rightarrow +\infty$  arises terms on the boundary  which are the same  than in~\cite{BL}, namely a distribution $\ell_1+\ell_2$ where $\ell_1$ is a measure absolutely continuous with respect to $\mu$ and   where, denoting by $N(x)$ the exterior  normal vector to $\partial\Omega$,  $\ell_2$  satisfies
$$\ell_2=  {\delta(\xi-\xi^+)-\delta(\xi-\xi^-)\over ( \xi_+-\xi_-)\cdot N(x)} \,\lambda\, {\bf 1}_{{\mathcal H}\cup{\mathcal G}}.$$
The reader can refer to~\cite{B} (pages 14-15) and~\cite{BL}  where the computations are carefully carried out.

On the other hand, if one uses the equation to transform $I_n$, there appears a new term which was not in the preceding articles. We are going to discuss why these terms are harmful when one has (3) of Assumptions~\ref{ass2}. This term involves the quantity $\nabla\cdot \left(\gamma(x) \theta_n\right)$. We introduce the principal symbol  $\tilde \gamma(y)\cdot \eta$  of the operator~$\Gamma$ which arises when writing $\nabla\cdot\left(\gamma\cdot\right)$ in the normal geodesic coordinates and we obtain
$$
I_n(t)  =  \left(q_b(y,D_{y'})u_n(t) \;,\;\tilde \gamma(y) \cdot \nabla \theta_n(t) \right) - \left( q_b(y,D_{y'}) \tilde \gamma(y) \cdot \nabla \theta_n(t)\;,\; u_n(t)\right) +o(1).
$$
Let us focus on the first term of $I_n$.
In order to study its contribution on the boundary, we use $\delta>0$ and a function $\chi\in{\mathcal C}^\infty_0(\R)$ such that
$\chi(t)=1$ for $t<1/2$ and $\chi(t) =0$ for $t>1$ with $0\leq \chi\leq 1$ and we study
$$J_{n,\delta}(t)= \left(\chi\left({y_d\over \delta}\right)q_b(y,D_{y'})u_n \;,\;\tilde \gamma(y) \cdot \nabla \theta_n \right) .$$
The fact that $\gamma$ is tangent to the boundary implies that
$$\tilde \gamma(y',0)=\left(\tilde \gamma'(y',0),0\right);$$
therefore,  we can write $\tilde \gamma=(\tilde \gamma', y_d\tilde\gamma_d)$. Then,
the worst term to estimate -- which is the one which involves $\partial_{y_d}$ derivatives of $\theta_n$ -- writes
$$\tilde J_{n,\delta}(t)= - \left( \chi\left({y_d\over \delta}\right)q_b(y,D_{y'})y_d \tilde \gamma_d(y) \partial_{y_d} \theta_n(t)\;,\; u_n(t)\right)$$
and we observe
\begin{eqnarray*}
\left|\tilde J_{n,\delta}(t)\right| & = & \left|\int_0^{+\infty}   \chi\left({y_d\over \delta}\right)\left(q_b(y,D_{y'})y_d \tilde \gamma_d(y) \partial{y_d} \theta_n(t)\;,\; u_n(t)\right)_{L^2(\R^{d-1}_{y'})} dy_d\right|\\
& \leq & \delta \|u_n(t)\|_{H^1(\Omega)} \| \partial_y\theta_n(t)\|_{L^2(\Omega)}.
\end{eqnarray*}
Therefore, for any $\Theta\in{\mathcal C}_0^\infty(\R)$, there exists a constant $C$
$$\left| \limsup_{n},{+\infty} \int \Theta(t) J_{n,\delta}(t) dt \right|\leq C\, \delta.$$
Finally, letting $\delta$ go to $0$, we obtain  that this term has no contribution on the boundary.
As a conclusion, we obtain $\ell=\ell_1+\ell_2$ and we can conclude like in~\cite{BL}.

\subsection{The link between the  measure $\mu$ and the measures $\mu_t^\pm$}\label{subsec:muetmupm}

In this section, we prove Proposition~\ref{muetmupm}.\\
\ni Note first that the second assertion is a simple consequence of the first one. Suppose that we have~\aref{eq:decompmu0} where,
by Proposition~\ref{prop:measure}, the measures $\mu^\pm_t$ satisfy
$$\partial_t \mu^\pm_t \pm \xi \cdot \nabla _x \mu^\pm_t =2 \alpha_\infty \mu^\pm _t,\;\;{\rm in}\;\;\Omega_1\times{\bf S}^{d-1}.$$
We obtain
\begin{eqnarray*}
H_{p_0}\mu & = & 2\left( \tau \partial_t \mu -\xi\cdot \nabla _x \mu\right)\\
& = &  \left( - \partial_t \mu_t^+ -\xi\cdot \nabla _x \mu_t^+\right )\otimes \delta(\tau+1)\otimes \d t\\
&
&\qquad
+ \left(  \partial_t \mu_t^- -\xi\cdot \nabla _x \mu_t^-\right)\otimes \delta(\tau-1)\otimes \d t\\
& = &- 2\alpha_\infty \mu^+_t(x,\xi) \otimes \delta(\tau+1)\otimes \d t + 2\alpha_\infty \mu^-_t(x,\xi) \otimes \delta(\tau-1)\otimes \d t \\
& = & 4\tau\, \alpha_\infty\, \mu.
\end{eqnarray*}
Let us now prove~\aref{eq:decompmu0}. We take $q\in {\mathcal A}_i^0$ and apply $q(z,D_z)$ to the first equation of~\aref{eq:thermo}. We get
$$0=\left(q(z,D_z) \left( \partial_t^2 u_n -\Delta u_n+\nabla\cdot(\gamma(x)\theta_n)\right)  ,u_n\right).$$
Using that $(\theta_n)_{n\in\N}$ is u.b. in $L^2_{loc}\left(\R, H^1(\Omega)\right)$ and that $u_n$ goes to $0$ weakly in $H^1(\Omega)$, we obtain passing to the limit
$$\int_{\R\times\Omega_1\times S^d(2)} q(z,\zeta)(|\xi|^2-\tau^2) d\mu(z,\zeta)=0.$$
Therefore, on the support of $\mu$, we have $\tau^2=|\xi|^2$. Since moreover $$|\zeta|^2=\tau^2+|\xi|^2=2\;\;{\rm on}\;\; S^d(2),$$ we obtain that $\tau^2=|\xi|^2=1$ on the support of $\mu$, whence the existence of two measures $\tilde \mu_\pm$ on $\R\times \Omega_1\times S^{d-1}$ such that
\begin{equation}\label{eq:decompmu}
\mu(t,x,\tau,\xi)=\delta(\tau+1)\otimes \tilde\mu_+(t,x,\xi) +\delta(\tau-1)\otimes \tilde\mu_-(t,x,\xi)\;\;{\rm in}\;\;{\mathcal D}'\left(\R\times\Omega_1\times S^d(2)\right).
\end{equation}
Let us now  link the measure $\tilde\mu_\pm $  with the measures~$\mu^\pm$ of Section~\ref{subsec:theonondeg} inside $\Omega$. Let us consider as in Section~\ref{subsec:theonondeg}
$$f_{\pm,n}(t,x)={1\over \sqrt 2}\left( i\partial_t u_n(t,x) \pm |D_x| u_n(t,x)\right).$$
These families are uniformly bounded in $L^2(\Omega)$ for all $t\in \R$ and their microlocal defect measures are the measures $\mu^\pm$.   Besides, by the definition of $\mu$, for $q\in{\mathcal A}_i^0$,
we have
\begin{eqnarray*}
\left(q(z,D_z) f_{\pm,n},f_{\pm,n}\right)  & =   &{1\over 2} \left( (i\partial_t \pm|D|)q(z,D_z) (i\partial_t \pm |D|)u_n,u_n\right)\\
& \td_n,{+\infty} & {1\over 2}  \int_{M\times S^d(2)} q_\infty(z,\zeta) (-\tau\pm|\xi|)^2 d\mu(z,\zeta)\\
 & & =  2 \int_{\R\times \Omega\times S^{d-1}} q_\infty(t,x,\mp1,\xi) d\tilde \mu_\pm(t,x,\xi)
 \end{eqnarray*}
 where we have used~\aref{eq:decompmu} for the last equality.
Let us choose $q(z,\zeta)=\chi(t) b(x,\xi)$ with $b(x,\xi)$ symbol of order $0$ compactly supported in $\Omega$, we obtain in the one hand
$$
 \int_{\R\times\Omega\times S^{d-1}} q_\infty(t,x,\mp1,\xi)  d\tilde \mu_\pm(t,x,\xi)
 = \int_{\R\times\Omega\times S^{d-1}} \chi(t)b(x,\xi)  d\tilde \mu_\pm(t,x,\xi).
$$
In the other hand, we have for all $t\in\R$,
$$(b(x,D)f_{\pm,n}(t),f_{\pm,n})\td_n,{+\infty}\int_{\Omega\times S^{d-1}} q(x,\omega) \d\mu^\pm_t(x,\omega).$$
Whence $\displaystyle{\tilde\mu_\pm(t,x,\omega)={1\over 2} \mu^\pm_t(x,\omega)\otimes \d t.}$

%%%%%%%%%%%%%%%%%%%%%%%%%%%%%%%%%%%%%%%%%%%%%%%%%

\end{document}